# TRANSVERSE STEADY BIFURCATION OF VISCOUS SHOCK SOLUTIONS OF A HYPERBOLIC-PARABOLIC MODEL IN A STRIP

RAFAEL A. MONTEIRO

ABSTRACT. In this article we derive rigorously a nonlinear, steady, bifurcation through spectral bifurcation (i.e., eigenvalues of the linearized equation crossing the imaginary axis) for a class of hyperbolic-parabolic model in a strip. This is related to "cellular instabilities" occuring in detonation and MHD. Our results extend to multiple dimensions the results of [AS12] on 1D steady bifurcation of viscous shock profiles; en passant, changing to an appropriate moving coordinate frame, we recover and somewhat sharpen results of [TZ08a] on transverse Hopf bifurcation, showing that the bifurcating time-periodic solution is in fact a spatially periodic traveling wave. Our technique consists of a Lyapunov-Schmidt type of reduction, which prepares the equations for the application of other bifurcation techniques. For the reduction in transverse modes, a general Fredholm Alternative-type result is derived, allowing us to overcome the unboundedness of the domain and the lack of compact embeddings; this result apply to general closed operators.

## Contents

## 1. Introduction

In this article we derive rigorously a nonlinear bifurcation through spectral bifurcation (i.e., eigenvalues in the linearized equation crossing the imaginary axis) for a multi-dimensional hyperbolic-parabolic model on a strip. We study a quasilinear, second order reaction diffusion equation that describes a fluid flowing in a duct with periodic boundary conditions in the y-direction. The equation is given by

$$u_t(x,y,t) + div f[\epsilon, u(x,y,t)] = \triangle u(x,y,t), \tag{1.1}$$

where $(x,y,t) \in \mathbb{R} \times \mathbb{T} \times \mathbb{R}$, $u : \mathbb{R} \times \mathbb{T} \times \mathbb{R} \to \mathbb{R}^n$, $f \in \mathscr{C}^2(\mathcal{I} \times \mathbb{R}^n; \mathbb{R}^n)$, smooth functions with all derivatives bounded, $\mathbb{T} = [0, 2\pi]$, $\mathcal{I}$ is an interval and $\epsilon$ the bifurcation parameter; in what follows, we denote $\mathbb{R} \times \mathbb{T}$ by $\Omega$. (1.1) can also can be rewritten as $u_t + [f^{(1)}(\epsilon, u)]_x + [f^{(2)}(\epsilon, u)]_y = \partial_x^2 u + \partial_y^2 u$.

1.1. **The mathematical setting.** Consider a family of viscous shock waves $\bar{u}^\epsilon$ solving (1.1), a traveling wave with speed c and y independent, i.e. $\bar{u}^\epsilon(x, y, t) = \bar{u}^\epsilon(\xi)$, where $\xi = x - ct$ and so that $\lim_{\xi \to \pm\infty} \bar{u}^\epsilon(\xi) = u^\pm(\epsilon)$. After a Galilean change of frame $(x, y, t) \mapsto (\xi, y, t)$ we can rewrite (1.1) as

$$v_t + [f^{(1)}(\epsilon, v) - cv]_x + [f^{(2)}(\epsilon, v)]_y - \triangle v = 0 \tag{1.2}$$

Assuming that $v + \bar{u}^\epsilon$ is a solution to (1.1) we derive a linearization about $\bar{u}^\epsilon$

$$\frac{\partial v}{\partial t} + [(Df^{(1)}(\epsilon, \bar{u}^\epsilon) - cI)v]_x + [Df^{(2)}(\epsilon, \bar{u}^\epsilon)v]_y - \triangle v = [A(\epsilon, \bar{u}^\epsilon)v]_x + [B(\epsilon, \bar{u}^\epsilon)v]_y - \triangle v = \mathscr{N}[\epsilon, \bar{u}^\epsilon, v],$$

where

$$\mathscr{N}(\epsilon, \bar{u}^\epsilon, v) = -[\mathscr{R}^{(1)}(\epsilon, \bar{u}^\epsilon, v)]_x - [\mathscr{R}^{(2)}(\epsilon, \bar{u}^\epsilon, v)]_y, \tag{1.3a}$$

$$\mathscr{R}^{(i)}(\epsilon, \bar{u}^\epsilon, v) = [f^{(i)}(\epsilon, v + \bar{u}^\epsilon) - f^{(i)}(\epsilon, \bar{u}^\epsilon) - Df^{(i)}(\epsilon, \bar{u}^\epsilon)v], \qquad i \in \{1, 2\} \tag{1.3b}$$

$(Df^{(1)}(\epsilon, \bar{u}^\epsilon) - cI)v = A(\epsilon, \bar{u}^\epsilon)v$ and $Df^{(2)}(\epsilon, \bar{u}^\epsilon)v = B(\epsilon, \bar{u}^\epsilon)v$. We can rewrite the above problem as

$$\frac{\partial v}{\partial t} + \mathscr{L}[\epsilon, \bar{u}^\epsilon, v] = \mathscr{N}[\bar{u}^\epsilon, w], \tag{1.4}$$

for $\mathscr{L}[\epsilon, \bar{u}^\epsilon, v] = [(Df^{(1)}(\bar{u}^\epsilon) - cI)v]_x + [Df^{(2)}(\bar{u}^\epsilon)v]_y$. It is clear that any y-translation of a solution to (1.4) is also a solution to it, a fact that will be used throughout the paper.

**Observation 1.** *[Translational symmetry] Let $v(\cdot, \cdot)$ be a solution of (1.4). Then $v(\cdot, \cdot + c)$ is also a solution of (1.4) for all $c \in \mathbb{R}$. In this case we say that "the system has translational symmetry in y-direction".*

As is common in some physical systems (see [Kno96]), one can study other kinds of symmetries besides the translational one; we say that solution u of (1.4) has rotational symmetry if there exists a constant matrix $R \in O(n)$ such that, if u solves (1.4) then $\Gamma[u]$ also solves it, where $\Gamma[u]$ is given by

$$\Gamma[u](x, y) = R \cdot u(x, -y). \tag{1.5}$$

**Definition 1.** *[O(2) symmetry] Consider a system of equations of the form (1.4) with traveling wave $\bar{u}^\epsilon$ that features translational symmetry in y-direction. If every solution is rotational symmetry under a fixed operator $\Gamma$ as defined in (1.5), if the operator commutes with $\Gamma$ upon linearization and if $\bar{u}^\epsilon$ is invariant under $\Gamma$ (i.e., $\Gamma[\bar{u}] = \bar{u}^\epsilon$) then the system (1.4) is said to have O(2) symmetry.*

As consequences of this definition, $v$ is a solution of (1.4) if and only if $\Gamma[v]$ is also a solution to the same equation; in other words:

**Observation 2.** *The linearization (1.4) preserves O(2) symmetry.*

Furthermore, $\lambda\phi$ are eigenvalue/eigenfunction of $\mathscr{L}$ if and only if $\lambda\Gamma(\phi)$ is an eigenvalue/eigenfunction of $\mathscr{L}$. In this article, we will be concerned with both cases:

    i. The system in (1.4) has translational symmetry in y-direction;
    ii. The system in (1.4) is $O(2)$ symmetric (under which the traveling wave is invariant).

Let's consider now a Galilean change of frame in y-direction $((\xi, y, t) \mapsto (\xi, y - dt, t))$ in order to make a perturbation $w$ in system (1.4) to be steady in the $O(2)$ case:

$$\frac{\partial v}{\partial t} + [(Df^{(1)}(\bar{u}^\epsilon) - cI)v]_x + [(Df^{(2)}(\bar{u}^\epsilon) - dI)v]_y - \triangle v = \mathscr{N}[\bar{u}^\epsilon, v], \tag{1.6}$$

Intuitively, if we reflect an object moving upward with velocity +d in the same axis its mirror image would be going downward with speed -d. So, if $d \neq 0$ there would be no way to "track" both objects ("real" + "mirror" image) at the same time in such a way that both of them would look steady; this is the main idea behind the next proof; now we make these words more rigorous.

**Observation 3.** *Reflection symmetry of system 1.4 is respected upon change of frame in y-direction if and only if d = 0.*

*Proof.* We know that as both $v(x, y)$ and $\Gamma[v] = Rv(x, -y)$ solve (1.2), therefore both $v(x - ct, y - dt, t)$ and $Rv(x - ct, -y - dt, t)$ solve (1.6). However, $\Gamma[v(x - ct, y - dt, t)]$ does not satisfy the equation in the moving frame, but instead, the equation

$$\frac{\partial w}{\partial t} + [(Df^{(1)}(\bar{u}^\epsilon) - cI)w]_x + [(Df^{(2)}(\bar{u}^\epsilon) + dI)w]_y - \triangle w = \mathcal{N}[\bar{u}^\epsilon, w],$$

Therefore, if we want both $v$ and $\Gamma[v]$ to solve the linearized equation we need d=0. □

**Remark 1.** *To set the dependence on the Galilean frame in y-direction explicitly, we will denote*

$$\mathscr{L}^d[w] := [(Df^{(1)}(\bar{u}^\epsilon) - cI)w]_x + [(Df^{(2)}(\bar{u}^\epsilon) - dI)w]_y - \triangle w \tag{1.7}$$

*minding that $d = 0$ whenever the system is $O(2)$ symmetric, in which case we will represent the operator by $\mathscr{L}^0$ or simply $\mathscr{L}$.*

We are looking for steady bifurcations of the system

$$\mathscr{L}^d[\epsilon, \bar{u}^\epsilon, w] = \mathcal{N}[\epsilon, \bar{u}^\epsilon, w] \tag{1.8}$$

1.2. **Plane waves.** A common approach to this problem in physics (as can be seen in [Cha61]) is a separation in the x-y dynamics in equation (1.8); one could represent a solution then as $w_k(x, y) = e^{iky}u_k(x)$. The $2\pi$ periodicity of solutions with respect to y implies that the wave numbers k's are in a discrete set, namely, $k \in \mathbb{Z}$. Let's omit the $\epsilon$ dependence for now; plugging in the linear part of (1.8),

$$\mathscr{L}^d[w_k] = e^{iky}\mathscr{L}^d_k[u_k] := e^{iky}\left((Au_k)' + ikBu_k + (k^2 - ikd)u_k - w_k''\right).$$

Here, $(\cdot)'$ denotes $\dfrac{d}{dx}$. Intuitively, we could take advantage of the boundedness of the domain in the y-direction to use Fourier analysis on it, hopefully recovering a solution to problem (1.8) by some sort of spectral synthesis; in fact, a substantial part of these notes is devoted to making this idea rigorous, specially the remainder of this section and **section 3**. In order to do that, we define the operators

$$\tilde{\Pi}_k[u](x, y) = \left(\frac{1}{2\pi}\int_0^{2\pi} u(x, \xi)e^{-ik\xi}\right) = \Pi_k[u](x)e^{iky}. \tag{1.9}$$

In fact, one can prove that $\tilde{\Pi}_k \circ \mathscr{L}[u](x, y) = \mathscr{L}^d_k \circ \Pi_k[u](x)e^{ik \cdot y}$ a.e., as shown in **appendix B**. Applying the operator $\tilde{\Pi}_k$ on both sides of equation (1.6) we get

$$\tilde{\Pi}_k \circ \mathscr{L}^d[w] = \mathscr{L}^d_k \circ \underbrace{\tilde{\Pi}_k[w]}_{w_k(x)e^{iky}} \underset{(1.3b)}{=} e^{iky}((Au_k)' + ikBu_k + (k^2 - ikd)_uk - u_k'') =$$

$$= \tilde{\Pi}_k \circ \mathcal{N}[u] = e^{iky}\underbrace{\left(\frac{1}{2\pi}\int_0^{2\pi} \mathcal{N}[u(x), \xi)]e^{-ik\xi}\right)}_{=: \mathcal{N}_k[u]} d\xi.$$

We end up with an infinite family of coupled 1-D problems:

$$\mathscr{L}_k^d[u](x) = \Pi_k \circ \mathscr{N}[u](x) = \mathscr{N}_k[u](x), \quad \text{for all} \quad k \in \mathbb{Z}, x \in \mathbb{R}. \tag{1.10}$$

Equation (1.11) is the essence of the plane waves method. We do not expect $\mathscr{N}$ and $\tilde{\Pi}_k$ to behave like $\mathscr{L}$ and $\tilde{\Pi}_k$, since $\mathscr{N}$ is nonlinear. The plane waves method allied to the Fourier synthesis will be in the backdrop of this work.

### 1.3. The bifurcation assumptions.
We will assume that the operators $\mathscr{L}_k^d(\epsilon, \bar{u}^\epsilon, \cdot)$ are all invertible for all $k \in \mathbb{Z} \setminus \{0, \pm k^*\}$, $k^* \neq 0$. We impose similar assumptions as in [TZ08a]:

(H0) There exists a simple eigenvalue $\lambda_{k^*}(\epsilon)$ (resp $\lambda_{-k^*}(\epsilon)$ ) associated to the point spectrum of $\mathscr{L}_{k^*}^d(\epsilon, \bar{u}^\epsilon, \cdot)$ (resp $\mathscr{L}_{-k^*}^d(\epsilon, \bar{u}^\epsilon, \cdot)$ ) such that, as $\epsilon \to 0$, $\mathscr{L}_{k^*}^d(\epsilon, \bar{u}^\epsilon, \cdot)$ (resp $\mathscr{L}_{-k^*}^d(\epsilon, \bar{u}^\epsilon, \cdot)$) goes through a spectral bifurcation, i.e., considering a bifurcation parameter $\epsilon \in \mathcal{I}$ and the mapping $\epsilon \mapsto \lambda_{\pm k^*}(\epsilon)$, we have that

$$Re(\lambda_{k^*}(0)) = 0 \quad , \quad Re(\lambda'_{k^*}(0)) \neq 0;$$

(H1) $(x, y) \mapsto f^{(j)}(x, y)$ is $\mathscr{C}^k$ on both entries, for $j \in \{1, 2\}$, $k \geq 2$; furthermore, in terms of $\mathscr{L}^2$ spectrum $\sigma(\cdot)$, there exists a $\beta < 0$ such that $Re(\sigma(\mathscr{L}_{\pm k^*}^{\bar{d}}[\epsilon]) \setminus \{\lambda_{\pm k^*}\}) < \beta$, for all $\epsilon \in \mathcal{I}$, where $\bar{d} = Im(\lambda_{k^*})$ (see **section 1.4** for the reasoning on the choice of $\bar{d}$);

(H2) Let $A(\epsilon, x) := df^1[\epsilon, \bar{u}^\epsilon(x)]$ and $A^\pm(\epsilon) = \lim_{x \to \pm\infty} A(\epsilon, x)$. We assume that $\sigma(A\pm(\epsilon))$ is real, distinct and nonzero. We denote the eigenspaces associated eigenvalues with positive (resp negative) real part by $\mathbb{E}^s(A^\pm)(\epsilon)$ (resp $\mathbb{E}^u(A^\pm)$). We assume that this system satisfy a Lax shock condition:

$$\dim \mathbb{E}^s(A^+(\epsilon)) \quad + \quad \dim \mathbb{E}^u(A^-(\epsilon)) \quad = \quad n+1, \quad \epsilon \in \mathcal{I}.$$

**Remark 2.** *As the endstate matrices are hyperbolic, (H2) clearly implies that $\bar{u}^\epsilon + \infty) \neq \bar{u}^\epsilon - \infty)$. Further, as described in [TZ08a, section 3.1 and lemma 3.2], differentible dependence on parameters gives that $(\epsilon, x) \mapsto \bar{u}^\epsilon(x)$ is $\mathscr{C}^2$ on both variables. Further, for $l \in \{0, 1, 2\}$, $|\partial_x^l(\bar{u}^\epsilon - \bar{u}^\epsilon \pm \infty))(x)| \leq Ce^{-\alpha|x|}$ as $x \to \pm\infty$ for some $\alpha > 0$. By the dependence of the matrices $A(\epsilon, \cdot)$ on $\bar{u}^\epsilon$ we have that $\partial_\epsilon^j A(\epsilon, x) \to \partial_\epsilon^j A^\pm(\epsilon)$ uniformly over compact subsets of $\mathcal{I}$ for every $j \in \{1, 2, \ldots k\}$.*

**Remark 3.** *On assumption (H0), the differentiability of the eigenvalues on $\epsilon$ is not asserted, but rather derived from the analytic dependence of the eigenvalues on the operators (in the operator norm, as explained in appendix C); therefore, this assumption is legitimate.*

**Remark 4.** *As the matrices $A^\pm(\cdot)$ are hyperbolic their eigenvalues depend continuously on the parameter $\epsilon$. we can assume then, without loss of generality, that the condition (H2) holds in such a way that both $\mathcal{I} \ni \epsilon \mapsto \dim(A^\pm(\epsilon)$ are constants.*

For our convenience in the simultaneous treatment of both problems, y-translational symmetric and O(2)-symmetric, we write the system (1.10) as

$$\mathscr{L}_k^{\bar{d}}[u](x) = \mathscr{N}_k[u](x) + ik(\bar{d} - d)u, \quad \text{for all} \quad k \in \mathbb{Z}, x \in \mathbb{R}, \tag{1.11}$$

minding that the last term on the right had side vanishes in the O(2)-symmetric case.

### 1.4. Escaping from the Hopf bifurcation scenario.
As we have seen in **section 1.3**, the spectral assumptions (H0) do not introduce constraints in the imaginary values of the crossing eigenvalues, abandoning the "Hopf-type bifurcation" scenario treated, for example, in [TZ08a, TZ05, TZ08b], which would assume $Im\{\lambda_{k^*}(0)\} \neq 0$. The main idea in order to avoid such a constraint lies in the Galilean change of reference frame in y-direction: one can choose a frame of reference in (1.6) such that the bifurcating eigenvalues happen to cross the imaginary axis at the origin. In fact, in the general, y-translational case, we know that if $(v(x), \lambda)$ are eigenfunction, eigenvalue of the system $\mathscr{L}_k^d[\cdot]$ then $(e^{iky+\lambda t}v(x), \lambda)$ is an eigenfunction of $\mathscr{L}_k^d[\cdot]$. What happens with the bifurcating eigenvalues if we change the frame of reference in y to $y - dt$? First of all, notice that $\mathscr{L}_k^d[w] := \mathscr{L}_k^0[w] - ikdw$. Clearly, if w,$\lambda$ are, respectively, eigenfunction,eigenvalue to $\mathscr{L}_k^0[w]$ then $\lambda - ikd$, w is and eigenfunction of $\mathscr{L}_k^d[w]$. Therefore, if we know that $\lambda - ikd$ is in the point spectrum of $\mathscr{L}_k^d$, the change of frame maps this problem to another one in which the eigenvalue is real in the system $\mathscr{L}_k^d[w]$. We could then map a Hopf bifurcation condition of a pair of eigenvalues/eigenfunctions $(\lambda \pm i\mu, v_\pm)$ crossing the imaginary axis at $\mu \neq 0$ to another problem in which there is a pair of eigenvalues crossing the imaginary axis through the origin Thus, although apparently less generic, we are actually dealing with a more general case than a Hopf bifurcation; overall, we could solve both cases at the same time.

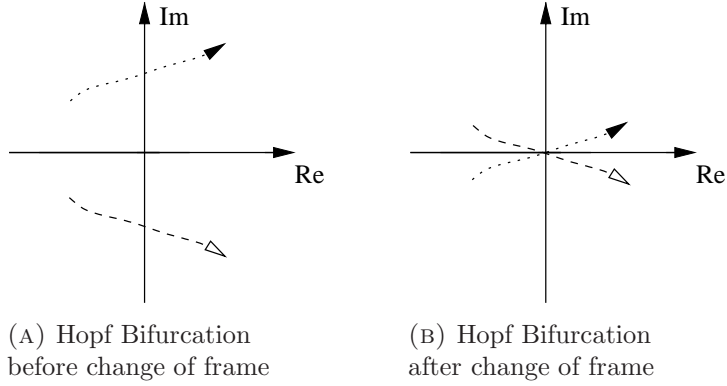

(A) Hopf Bifurcation before change of frame

(B) Hopf Bifurcation after change of frame

FIGURE 1. Vertical translation of the crossing eigenvalues after a Galilean change of frame in y-direction

**Remark 5.** *Notice that the Galilean change of frame in y-direction, when applied to a standard Hopf-bifurcation, introduces a $d \neq 0$ in equation (1.7). Therefore, in virtue of the observation 3, the technique presented in section 1.4 breaks O(2) symmetry.*

1.5. **Results.** Assumption (H0) does not impose any constraint in the imaginary part of the eigenvalues $\lambda_{\pm k^*}$. Taking into account the reasoning in **section 1.4**, it is shown in section 4 that we can choose a change of frame in y-direction so that bifurcating eigenvalues occur through the origin; namely, we choose a function $\epsilon \mapsto d(\epsilon)$ in equation (1.6) such that $d(0) = Im(\lambda_{k^*}(0))$. As explained above, $\mathscr{L}_{\pm k^*}^{d(0)}$ has bifurcating eigenvalues through the origin only. In **section 6** we prove the following:

**Theorem 1.12.** *[Translation invariant case] Let system (1.11) satisfying the bifurcation assumptions (H0-H2) have translation symmetry. Then, there exists an interval $\mathcal{J} \subset \mathbb{R}$ so that $0 \in \mathcal{J}$, a mapping $\mathcal{J} \ni s \mapsto v(s) \in \mathscr{H}^2(\Omega; \mathbb{R}^n)$ and a continuous mapping $\mathcal{J} \ni s \mapsto d(s) \in \mathbb{R}$ such that*

$$\mathcal{J} \ni s \mapsto \mathscr{L}^{d(s)}[v(s)] = \mathscr{N}[\bar{u}^\epsilon, v(s)],$$

*where $v(s) \not\equiv 0$, $v(0) = 0$ and $d(0) = \bar{d}$, for all $|d(s) - \bar{d}|$ sufficiently small in a neighborhood of s=0*

As a corollary of this result, in the "Hopf-bifurcation scenario" $(d(0) = Im(\lambda_{k^*}) \neq 0)$, we could obtain periodic solutions in the y-direction, partially recovering results of [PYZ03]:[1]

---

[1]In [PYZ03] there is determined the entire bifurcation diagram, including not only the traveling waves demonstrated here, but also reflectionally symmetric time-periodic solutions given by a nonlinear superposition of counter-propagating traveling waves.

**Corollary 1.** *[O(2) Hopf-bifurcation] Let system* (1.11) *satisfying the bifurcation assumptions (H0-H2) have O(2) symmetry. Assume also assume that the eigenvalues $\lambda_{\pm k^*}(0)$ are double [2] and $Im(\lambda_{\pm k^*}(0)) \neq 0$. Then, there exists an interval $\mathcal{J} \subset \mathbb{R}$, $0 \in \mathcal{J}$, a mapping $\mathcal{J} \ni s \mapsto v(s) \in \mathscr{H}^2(\Omega;\mathbb{R}^n)$ and a continuous mapping $\mathcal{J} \ni s \mapsto d(s)$ such that*

$$\mathcal{J} \ni s \mapsto \mathscr{L}^{d(s)}[v(s)] = \mathscr{N}[\bar{u}^\epsilon, v(s)],$$

*where $v(s) \not\equiv 0$, $v(0) = 0$ and $d(0) = \bar{d}$, for all $|d(s) - \bar{d}|$ sufficiently small in a neighborhood of s=0. Furthermore, v is a periodic function in y-direction with period $\tau(s) = \dfrac{2\pi}{d(s)}$.*

Analogously, we prove

**Theorem 1.13.** *[O(2) symmetric case] Let system* (1.11) *satisfying the bifurcation assumptions (H0-H2) have also O(2) symmetry. Then, there exists an interval $\mathcal{J} \subset \mathbb{R}$ so that $0 \in \mathcal{J}$, a mapping $\mathcal{J} \ni s \mapsto v(s) \in \mathscr{H}^2(\Omega;\mathbb{R}^n)$ and a continuous mapping $\mathcal{J} \ni s \mapsto d(s)$ such that*

$$\mathcal{J} \ni s \mapsto \mathscr{L}^{d(s)}[v(s)] = \mathscr{N}[\bar{u}^\epsilon, v(s)],$$

*where $v(s) \not\equiv 0$, $v(0) = 0$ and $d(0) = \bar{d}$, for all $|d(s) - \bar{d}|$ sufficiently small in a neighborhood of s=0.*

In passing, a Fredholm Alternative-type theorem is derived in the **appendix D**, allowing us to overcome the unboundedness of the domain and the lack of compact embeddings (here, in the case $\mathscr{H}^2(\mathbb{R};\mathbb{R}^n) \hookrightarrow \mathscr{L}^2(\mathbb{R};\mathbb{R}^n)$); this part is crucial to reduce the problem to a finite dimensional setting. Furthermore, we make use of plane waves as suggested above in a rigorous manner, encoding not only a separation of dynamics/variables technique, but rather a spectral synthesis in the periodic direction.

1.6. **Discussion and open problems.** The study of steady longitudinal (or 1D) bifurcations of shock waves was carried out in [AS12]. An ODE problem involving bifurcation of heteroclinic orbits connecting hyperbolic rest points, this was treated by introducing a Melnikov separation function between stable and unstable manifolds and studying an associated mapping bifurcation problem for its zeros. Our analysis can be viewed as the natural generalization to steady transversal (or multi-D) bifurcation.

It is interesting to note that our treatment of the $k = 0$ modes corresponding to longitudinal perturbations proceeds by essentially the same steps as in [AS12], but in a different order. Namely, in the ODE case, one integrates the traveling-wave ODE to obtain a first-order system with hyperbolic rest points, introduces a Melnikov separation function by taking a section, then studies the resulting bifurcation problem by linearization and Lyapunov–Schmidt reduction; in our case, on the other hand, we first linearize, integrate to a first-order system with hyperbolic behavior at infinity, then we introduce a phase condition effectively taking a section, after which we finally perform a Lyapunov–Schmidt reduction. This connection, noted in hindsight, may provide useful orientation for the analysis of longitudinal modes $k = 0$. The analysis may equally well be motivated as a "bordering method" in the study of general boundary-value problems with nonzero Fredholm index, recovering under additional "phase conditions" resolvent estimates similar to those of the Fredholm zero case (as described, for example, in [Hen81]).

Both types of bifurcation can occur in magnetohydrodynamics (MHD), the longitudinal type as described in [AS12] and the transverse type as described in Example 1 below.

Our results at the same time generalize the study of transverse Hopf-type bifurcation carried out in [TZ08a, section 6.1], in which bifurcating eigenvalues were assumed to have nonzero imaginary part, and a resulting time-periodic nonlinear bifurcation was established. As described above, the value of the imaginary parts may be prescribed arbitrarily by choice of an appropriate moving coordinate frame in the y-direction; hence in some sense the two problems are equivalent. However, the solution is obtained here by quite different and more direct techniques, which yield also the additional information that the time-periodic solutions of [TZ08a] are, more precisely, spatially-periodic traveling waves.

---

[2]Due to O(2) symmetry, the multiplicity of eigenvalues is always at least two; see discussion in section 6.3.

Our results here are obtained for "artificial," Laplacian viscosity. However, they readily generalize to quasilinear strictly parabolic systems. The solution to problem (1.1) with real viscosity is challenging and still open. Naturally, other types of cross sectional symmetric domains must provide very interesting study cases and different kinds of bifurcations. Nonlinear stability of bifurcating solutions is another interesting open problem.

**Example 1** (Transverse bifurcation in MHD). *In [FT08], Freistühler and Trakhinin exhibit an example of a "parallel" MHD shock wave (meaning one for which the magnetic field direction is parallel to the normal to the shock surface), which at the inviscid level has when considered on the whole space a transverse instability, i.e., a dispersion relation $\lambda(k)$ with $Re\lambda_k > 0$ for $k \neq 0$, but at the viscous level is one-dimensionally (i.e., for $k \neq 0$) stable. More, the eigenvalue $\lambda(k)$ is simpe for each $k \neq 0$, whence, by $O(2)$ symmetry, it must be real, and thus positive. But, by the results of [SZ99] relating viscous and inviscid spectra, this implies that there is a nearby positive real root $\lambda_*(k)$ that is tangent to $\lambda(k)$, i.e., a viscous transverse instability as well nonresonance. Considering the same shock as a solution on a strip of width $L$ with periodic boundary conditions, and treating $L$ as a control parameter, we find that only wave numbers $k = 2\pi j/L$, $j$ integer, are allowable. Thus, taking $L \to 0$, we have a larger and larger spectral gap between $k = 0$ and the next mode $k = 2\pi/L$. Since high-frequency estimates imply that all spectra are stable (negative real part) for $|k|$ sufficiently large, this means that the shock eventually becomes stable on the bounded cross-section domain, hence there must be some point at which an eigenvalue crosses the imaginary axis at $k = 2\pi/L$, a transverse bifurcation. This may be of steady, or Hopf type. To conclude nonlinear bifurcation, we need the further conditions if simplicity and nonresonance*

## 2. A few technical lemmas

### 2.1. On the notation used.
In this paper, we define $\mathbb{Z}^* := \mathbb{Z} \setminus \{0, \pm k^*\}$ and denote the domain $\mathbb{R} \times \mathbb{T}$ by $\Omega$. We use the following inner product for $u, v \in \mathscr{H}^s(X, \mathbb{R}^n)$:

$$\langle u, v \rangle_{\mathscr{H}^s(X,\mathbb{R}^n)} = \sum_{|\alpha| \leq s} \int_X \partial^\alpha[u(x)] \cdot \overline{\partial^\alpha[v(x)]} dx$$

Where "·" denotes the usual n dimensional inner product and $X$ is either $\mathbb{R}$ or $\Omega$. We frequently suppress the indication $\mathscr{H}^s(X, \mathbb{R}^n)$ when using the norm $||\cdot||_{\mathscr{H}^s(X;\mathbb{R}^n)}$, writing instead $||\cdot||_{\mathscr{H}^s(X)}$ to describe these objects. We intend then to keep the full meaning of the symbols without any risk of misunderstandings, avoiding unnecessary or cumbersome notation.

### 2.2. On the family of projections $\{\Pi_k[\cdot]\}_{k \in \mathbb{Z}}$.
To understand the relation between $w \in \mathscr{H}^2(\Omega; \mathbb{R}^n)$ and the family $\{w_k\}_{k \in \mathbb{Z}}$ of its projections we need to understand the operators $\{\Pi_k[\cdot]\}_{k \in \mathbb{Z}}$.

**Lemma 1.** *Let $\Pi_k(\cdot)$ be as defined in the previous session. Then we have that $\Pi_k : \mathscr{H}^s(\Omega; \mathbb{R}^n) \to \mathscr{H}^s(\mathbb{R}; \mathbb{C}^n)$. In addition to that, the following equality holds:*

$$||\partial_x^a \partial_y^b u||^2_{\mathscr{L}^2(\Omega)} = \sum_{k \in \mathbb{Z}} k^{2 \cdot b} ||\partial_x^a u_k||^2_{\mathscr{L}^2(\mathbb{R})},$$

*where $\partial_x u_k$ denotes the weak derivative of $u_k$ and $a, b \in \{0, 1 \ldots, s\}$, which we also show to exist.*

*Proof.* Proposition 16 on appendix B provides the existence of weak derivatives $u_k$; namely, $\partial^\alpha(\widetilde{\Pi}_k[u]) = (ik)^{\alpha_y}\widetilde{\Pi}_k[\partial^\alpha u]$; one can conclude that

$$||\partial_x^{\alpha_x}\partial_y^{\alpha_y}[u]||_{\mathscr{L}(\Omega;\mathbb{R}^n)} = \int_\Omega |\partial_x^{\alpha_x}\partial_y^{\alpha_y}[u]|^2 dxdy = \int_\mathbb{R}\left\{\int_\mathbb{T}|\partial_x^{\alpha_x}\partial_y^{\alpha_y}[u]|^2 dy\right\}dx =$$

$$\underset{Parseval's\ ineq.}{=} \int_\mathbb{R}\left\{\sum_{k\in\mathbb{Z}}\left|(ik)^{\alpha_y}\widetilde{\Pi}_k[\partial_x^{\alpha_x}u]\right|^2(x)\right\}dx = \int_\mathbb{R}\sum_{k\in\mathbb{Z}}|k|^{2\alpha_y}|\widetilde{\Pi}_k[\partial_x^{\alpha_x}u]|^2 dx =$$

$$\underset{Mon.conv.thm.}{=} \sum_{k\in\mathbb{Z}}|k|^{2\alpha_y}||\Pi_k[\partial_x^{\alpha_x}u]||^2_{\mathscr{L}(\mathbb{R})}.$$

□

**Lemma 1** suggests that defining a norm on the family of functions $(u_k)_{(k\in\mathbb{Z})}$ could be a good way to understand $u$ when $u \in \mathscr{H}^2(\mathbb{R}\times\mathbb{T})$; we will discuss this part more deeply in the **section 3**.

2.3. **The embedding** $\mathscr{H}^2(\Omega;\mathbb{R}^n) \subset \mathscr{L}^\infty(\Omega;\mathbb{R}^n)$. The solutions to (1.11) will be found through a contraction mapping theorem application. In order to do that, we will need some estimates in the nonlinear term (1.3a). The main result, which will be used throughout the paper, consists in proving that $u \in \mathscr{H}^2(\Omega;\mathbb{R}^n) \subset \mathscr{L}^\infty(\Omega;\mathbb{R}^n)$ continuously; the idea is a variation on the Sobolev lemma (see [Yos95, chapter VI, section 7]).

**Lemma 2.** *For almost every* $x \in \mathbb{R}$, $u(x,y) - u(x,z) = \int_z^y \partial_y u(x,w)dw$

*Proof.* Let $V(x,y) := \int_{y_0}^y \partial_y u(x,w)dw$; one can show that $V(x,\cdot)$ is finite.[3] We will show that the y-weak derivative of V(x,y) - u(x,y) is zero, so the later is y independent:

$$\int_\mathbb{R}\int_\mathbb{T} V(x,y)\cdot\partial_y\phi(x,y)dydx =$$

$$= \int_\mathbb{R}\left\{\int_0^{y_0}[\int_{y_0}^y \partial_y u(x,w)dw]\cdot\partial_y\phi(x,y)dy + \int_{y_0}^{2\pi}[\int_{y_0}^y \partial_y u(x,w)dw]\cdot\partial_y\phi(x,y)dy\right\}dx =$$

$$= \int_\mathbb{R}\left\{-\int_0^{y_0}[\int_y^{y_0}\partial_y u(x,w)dw]\cdot\partial_y\phi(x,y)dy + \int_{y_0}^{2\pi}[\int_{y_0}^y \cdot\partial_y u(x,w)dw]\cdot\partial_y\phi(x,y)dy\right\}dx =$$

$$\underset{reparametrizing}{=} \int_\mathbb{R}\left\{-\int_0^{y_0}[\int_0^w \partial_y u(x,w)\cdot\partial_y\phi(x,y)dy]dw + \int_{y_0}^{2\pi}[\int_w^{2\pi}\partial_y u(x,w)\cdot\partial_y\phi(x,y)dy]dw\right\}dx =$$

$$= \int_\mathbb{R}\left\{-\int_0^{y_0}\partial_y u(x,w)\cdot[\phi(x,w)-\phi(x,0)]]dw + \int_{y_0}^{2\pi}\partial_y u(x,w)\cdot[\phi(x,2\pi)-\phi(x,w)]dw\right\}dx =$$

$$\underset{\phi(x,0)=\phi(x,2\pi)}{=} -\int_\mathbb{R}\int_\mathbb{T}\partial_y u(x,w)\cdot\phi(x,w)dw + \underbrace{\int_\mathbb{R}\int_\mathbb{T}\partial_y u(x,w)\cdot\phi(x,0)dwdx}_{=0}$$

---

[3]Basically, it consists in proving that $\dfrac{\partial\hat{u}}{\partial y}(x,y) = \widehat{\dfrac{\partial u}{\partial y}}(x,y)$, in the weak sense, which mostly follows by the definition.

Therefore,
$$\int_\mathbb{R} \int_\mathbb{T} [V(x,y) - u(x,y)] \cdot \partial_y \phi(x,y) dy dx = 0 \Rightarrow V(x,y) - u(x,y) = C(x)$$
$$\Rightarrow \quad u(x,y) - u(x,z) = V(x,y) - V(x,z) = \int_z^y \partial_y u(x,w) dw.$$
□

**Theorem 2.1.** *The embedding $u \in \mathscr{H}^2(\Omega; \mathbb{R}^n) \subset \mathscr{L}^\infty(\Omega; \mathbb{R}^n)$ holds continuously.*

*Proof.* Denote $\frac{1}{2\pi} \int_\mathbb{T} u(x,w) dw$ by $M(x)$. Considering $u(x,y) = \{u(x,y) - M(x)\} + \{M(x)\} = I + II$ we'll show that $I$ and $II$ are bounded by a constant times $||u||_{\mathscr{H}^2(\Omega)}$.

$$I = u(x,y) - M(x) = \frac{1}{2\pi} \int_\mathbb{T} [u(x,y) - u(x,w)] dw \stackrel{lemma\ 2}{=} \frac{1}{2\pi} \int_\mathbb{T} \int_w^y \frac{\partial u}{\partial y}(x,z) dz dw =$$

$$\frac{1}{2\pi} \int_\mathbb{T} \int_\mathbb{T} 1_{[w,y]}(z) \frac{\partial u}{\partial y}(x,z) dz dw = \frac{1}{2\pi} \int_\mathbb{T} \int_\mathbb{T} 1_{[w,y]}(z) \left\{ \int_\mathbb{R} \widehat{\frac{\partial u}{\partial y}}(\cdot,z)(\xi) e^{ix\xi} d\xi \right\} dz dw$$

$$\Rightarrow |u(x,y) - M(x)| \lesssim ||\widehat{\frac{\partial u}{\partial y}}||_{\mathscr{L}(\Omega;\mathbb{R}^n)^1} \lesssim ||u||_{\mathscr{H}^2(\Omega)}.$$

where $1_\omega$ denotes the characteristic function of the set $\omega$. Moreover,
$$II = \frac{1}{2\pi} \int_\mathbb{T} u(x,w) dw = \frac{1}{2\pi} \int_\mathbb{T} \left\{ \int_\mathbb{R} \widehat{u(\cdot,w)}(\xi) e^{ix\xi} d\xi \right\} dw.$$

Thus, $|II| \lesssim ||u||_{\mathscr{H}^2(\Omega)}$, and therefore, $|u(x,y)| \lesssim ||u||_{\mathscr{H}^2(\Omega)} + |II| \lesssim ||u||_{\mathscr{H}^2(\Omega)}$. □

**2.4. A rate of decay in $|k|$ for the norm of the solution map $(\mathscr{L}_k^d)^{-1}$ as $|k| \to \infty$.** By the definition of $\mathscr{L}_k^d[\cdot]$ in (1.9) we derive through the weak formulation a bilinear form $B_k^d[\cdot,\cdot] : \mathscr{H}^1(\Omega; \mathbb{R}^n) \times \mathscr{H}^1(\Omega; \mathbb{R}^n) \to \mathbb{C}$

$$B_k^d[u,w] := \int_\mathbb{R} \left\{ -Au \cdot \bar{w}' + i(k+d) Bu \cdot \bar{w} + k^2 u \cdot \bar{w} + u' \cdot \bar{w}' \right\} dx = \int_\mathbb{R} f \cdot \bar{w},$$

$\forall w \in \mathscr{C}_c^\infty(\mathbb{R})$. In order to use Lax-Milgran's theorem we need to show coercivity and boundedness of this bilinear form[4]; we can prove then the following theorem:

**Theorem 2.2.** *$Re \circ B^d[\cdot,\cdot]$ is bounded and coercive for all $k$ sufficiently large. Moreover, the following estimates hold:*

$$||\mathscr{L}_k^{-1}(f)||_{\mathscr{L}^2(\mathbb{R}) \to \mathscr{L}^2(\mathbb{R})} = \mathcal{O}\left(\frac{||f||_{\mathscr{L}^2(\mathbb{R})}}{k^2}\right), \tag{2.3a}$$

$$||\partial_x \mathscr{L}_k^{-1}(f)||_{\mathscr{L}^2(\mathbb{R}) \to \mathscr{L}^2(\mathbb{R})} = \mathcal{O}\left(\frac{||f||_{\mathscr{L}^2(\mathbb{R})}}{k}\right), \tag{2.3b}$$

$$||\partial_x \mathscr{L}_k^{-1}(f)||_{\mathscr{H}^1(\mathbb{R}) \to \mathscr{L}^2(\mathbb{R})} = \mathcal{O}\left(\frac{||f||_{\mathscr{L}^2(\mathbb{R})}}{|k|} \wedge \frac{||f||_{\mathscr{H}^1(\mathbb{R})}}{k^2}\right), \tag{2.3c}$$

$$||\partial_x^2 \mathscr{L}_k^{-1}(f)||_{\mathscr{H}^1(\mathbb{R}) \to \mathscr{L}^2(\mathbb{R})} = \mathcal{O}\left(||f||_{\mathscr{L}^2(\mathbb{R})} \wedge (\frac{||f||_{\mathscr{H}^1(\mathbb{R})}}{|k|})\right). \tag{2.3d}$$

---
[4]Notice that if a bilinear form $B_k^d[\cdot,\cdot] : H \times H \to \mathbb{C}$ is complex valued, the Lax-Milgran theorem can be proved using the assumption that $Re \circ B_k^d[\cdot,\cdot]$ satisfies boundedness and coercivity; just use the fact that $Re \circ B[u,v] = Im \circ B[u,iv]$

*Proof.* As $\mathscr{L}_k^a[w] = \mathscr{L}_k^b[w] + ik(a-b)w$ for a, b real numbers, $B_k^d[w,v] = B_k^{\bar{d}}[w,v] + ik(d-\bar{d})\langle w,v\rangle$. Consequently, $Re\{B_k^d[w,w]\} = Re\{B_k^{\bar{d}}[w,w]\}$, i.e., $ik(d-\bar{d})$ has no role in the proof of invertibility of $\mathscr{L}_k^d[\cdot]$; in the proof we will fix $d = \bar{d}$. The boundedness of the bilinear form follows by successive aplication of Holder's inequality and other convexity inequalities, so we will just prove the coercivity. To short our notation, our norms will be denoted by $||\cdot||_{\mathscr{L}^2(\mathbb{R})} = ||\cdot||$ and $||\cdot||_{\mathscr{L}^\infty(\mathbb{R})} = ||\cdot||_\infty$:

$$Re(B_k^d[w,w]) \geq k^2||w||^2 + ||w'||^2 - \frac{(||A||\cdot||w||)^2}{2} - \frac{||w'||^2}{2} - k||B||\cdot||w||^2$$

While in the r.h.s. $Re\langle f,w\rangle \leq ||f||\cdot||w|| \leq \frac{||f||^2}{2\cdot k^2} + \frac{k^2||w||^2}{2}$. Putting these estimates together,

$$\left(\frac{k^2}{2} - k||B|| - \frac{||A||^2}{2}\right)||w||^2 + \frac{||w'||^2}{2} \leq \frac{||f||^2}{2\cdot k^2}$$

hence (2.3a) and (2.3b) follow. Analogously, assuming $f \in \mathscr{H}^1(\mathbb{R};\mathbb{R}^n)$, we have

$$\partial_x \mathscr{L}_k[w] = \mathscr{L}_k[\partial_x w] + (A'(\bar{u})\bar{u}'w')' + i(kB'(\bar{u})\bar{u}')w = f' \tag{2.4}$$

In the weak formulation, it results in $Re(B_k[w',w']) = Re(\langle f',w'\rangle + \langle A'(\bar{u})\bar{u}'w,w''\rangle - i\langle B'(\bar{u})\bar{u}'w,w'\rangle)$. Using Holder inequality and the estimates on remark 2 we obtain an upperbound for the left hand side for $\eta > 0$ sufficiently small:

$$Re(\langle f',w'\rangle) + ||A'(\bar{u})\bar{u}'||_\infty ||w||\cdot||w''|| + ||B'(\bar{u})\bar{u}'||_\infty ||w||\cdot||w'||$$

Rearranging, we obtain

$$\left((1-\eta)\frac{k^2}{2} - k||B|| - \frac{||A||^2}{2}\right)||w'||^2 + \frac{||w''||^2}{4} \leq Re(\langle f',w'\rangle) + C||w||^2$$

Now, we can estimate the term $Re(\langle f',w'\rangle)$ in two ways:

$$|Re(\langle f',w'\rangle)| \leq \frac{||f'||^2}{2\eta k^2} + \eta k^2 \frac{||w'||^2}{2} \tag{2.5a}$$

$$|Re(\langle f',w'\rangle)| \leq 2||f||^2 + \frac{||w''||^2}{8} \tag{2.5b}$$

In both cases, using the inequalities (2.3b) and (2.3c), we obtain:

$$\left((1-2\eta)\frac{k^2}{2} - k||B|| - \frac{||A||^2}{2}\right)||w'||^2 + \frac{||w''||^2}{4} \lesssim_\eta \frac{||f'||^2}{k^2} + ||w||^2 \lesssim_\eta \frac{||f||^2_{\mathscr{H}^1(\mathbb{R})}}{k^2}$$

$$\left((1-\eta)\frac{k^2}{2} - k||B|| - \frac{||A||^2}{2}\right)||w'||^2 + \frac{||w''||^2}{8} \lesssim_\eta ||f||^2 + ||w||^2 \lesssim_\eta ||f||^2_{\mathscr{L}^2(\mathbb{R})}$$

From the last inequalities we obtain estimates (2.3c) and (2.3d). $\square$

In fact, looking more carefully to the proof of the inequalities (2.3c) and (2.3d) one can prove the stronger result that, if $f \in \mathscr{L}^2(\mathbb{R};\mathbb{R}^n))$, then $||u_k||_{\mathscr{H}^2(\mathbb{R})} \lesssim ||f||_{\mathscr{L}^2(\mathbb{R})}$. Notice that it consists in showing that the assumption $f \in \mathscr{H}^1(\mathbb{R};\mathbb{R}^n)$ is not necessary.

**Proposition 1.** *Let $\mathscr{L}_k$ be an operator as defined in **section 1.2** with the property that $||u_k||_{\mathscr{H}^1(\mathbb{R})} \lesssim ||f||_{\mathscr{L}^2(\mathbb{R})}$[5]. Then $u \in \mathscr{H}^2(\mathbb{R};\mathbb{R}^n)$ and $||u_k||_{\mathscr{H}^2} \lesssim ||f||_{\mathscr{L}^2}$.*

---
[5]This condition is provided by equations (2.3a) and (2.3b) for large k and, for small k, by the results in **section 2.5**

*Proof.* For the sake of clarity we will omit the indexes for now. The idea consists in mollifying f and solving the problem $\mathscr{L}[u^\eta] = f^\eta$. As $\eta \to 0^+$, $f^\eta \to f$ in $\mathscr{L}^2(\mathbb{R}; \mathbb{R}^n)$. Showing that $u^\eta{}_{(\eta>0)} \subset \mathscr{H}^2(\mathbb{R}; \mathbb{R}^n)$ is Cauchy $u^\eta \to u \in \mathscr{H}^1(\mathbb{R}; \mathbb{R}^n)$ where $\mathscr{L}[u] = f$ implies that $u \in \mathscr{H}^2(\mathbb{R}; \mathbb{R}^n)$. Initially, notice that $||f^\eta||_{\mathscr{L}^2(\mathbb{R})} \leq ||f||_{\mathscr{L}^2(\mathbb{R})}$. (2.3c) and (2.3d) provide $||u^\eta||_{\mathscr{H}^2(\mathbb{R})} \lesssim ||f^\eta||_{\mathscr{L}^2(\mathbb{R})} \leq ||f||_{\mathscr{L}^2(\mathbb{R})}$. Thus $u^\eta{}_{(\eta>0)}$ is an $\mathscr{H}^2(\mathbb{R})$ bounded set. By linearity,

$$\mathscr{L}[u - u^\eta] = f - f^\eta \overset{||u_k||_{\mathscr{H}^1} \lesssim ||f||_{\mathscr{L}^2}}{\Rightarrow} ||u - u^\eta||_{\mathscr{H}^1(\mathbb{R})} \lesssim ||f - f^\eta||_{\mathscr{L}^2(\mathbb{R})} \overset{\eta \to 0}{\to} 0$$

and then $u^\eta \to u$ in $\mathscr{H}^1(\mathbb{R})$. Again by linearity we get

$$\mathscr{L}[u^\eta - u^\delta] = f^\eta - f^\delta \quad \Rightarrow \quad ||u^\eta - u^\delta||_{\mathscr{H}^2} \lesssim ||f^\eta - f^\delta||_{\mathscr{L}^2(\mathbb{R})} \to 0 \quad \text{as} \quad \eta, \delta \downarrow 0.$$

Therefore $u^\eta \to v$ in $\mathscr{H}^2(\mathbb{R}; \mathbb{R}^n)$, by completeness. In fact, $u = v$, because for any $|\alpha| \leq 2$ multi-index and $\phi \in \mathscr{C}_c^\infty(\mathbb{R}; \mathbb{R}^n)$,

$$(-1)^{|\alpha|} \int_\mathbb{R} u \partial^\alpha \phi \, dx = (-1)^{|\alpha|} \lim_{\eta \downarrow 0} \int_\mathbb{R} u^\eta \partial^\alpha \phi \, dx = \lim_{\eta \downarrow 0} \int_\mathbb{R} (\partial^\alpha u^\eta) \phi \, dx = \int_\mathbb{R} (\partial^\alpha v) \phi \, dx$$

□

2.5. **Estimates for small $|k|$.** Let $\mathscr{L}_k$ be an operator as defined in section 2.2 and $u_k$ a solution to the problem. $\mathscr{L}_k[u_k] = f$ in the weak sense. Suppose that $\mathscr{L}_k$ is invertible from $\mathscr{L}^2(\mathbb{R}; \mathbb{R}^n)$ to itself, i.e.,

$$||u_k||_{\mathscr{L}^2(\mathbb{R})} \lesssim ||f||_{\mathscr{L}^2(\mathbb{R})}. \tag{2.7}$$

These assumptions give possibility for an improvement of these estimate to higher order Sobolev spaces, as shown in the lemma below:

**Lemma 3.** *Let $u_k$ be a weak solution to $\mathscr{L}_k[u_k] = f$, where $||u_k||_{\mathscr{L}^2(\mathbb{R})} \lesssim ||f||_{\mathscr{L}^2(\mathbb{R})}$. Then $||u_k||_{\mathscr{H}^1(\mathbb{R})} \lesssim ||f||_{\mathscr{L}^2(\mathbb{R})}$. Furthermore*

$$||u_k||_{\mathscr{H}^2(\mathbb{R})} \lesssim ||f||_{\mathscr{L}^2(\mathbb{R})}. \tag{2.8}$$

*Proof.* The proof of theorem 2.2 shows that for every k the operator $(\mathscr{L}_k + \lambda I)[\cdot]$ is coercive in $\mathscr{H}^1(\mathbb{R}; \mathbb{R}^n)$ for $\lambda > 0$ sufficiently large. Notice then that $\mathscr{L}_k[u_k] = f \Leftrightarrow (\mathscr{L}_k + \lambda I)[u_k] = (f + \lambda u_k)$. Denoting by $B_\lambda[\cdot, \cdot]$ the bilinear form associated with the operator $(\mathscr{L}_k + \lambda I)$, we get ($c > 0$)

$$c||u_k||^2_{\mathscr{H}^1(\mathbb{R})} \leq \text{Re} \circ B_\lambda[u, u] \leq |\langle f + \lambda u_k, u_k \rangle| \leq ||f||_{\mathscr{L}^2(\mathbb{R})} \cdot ||u_k||_{\mathscr{L}^2(\mathbb{R})} + \lambda ||u_k||^2_{\mathscr{L}^2(\mathbb{R})}$$

Using the assumption (2.7), the right hand side can estimated by $||f||_{\mathscr{L}^2(\mathbb{R})}||u_k||_{\mathscr{L}^2(\mathbb{R})} + \lambda ||u_k||^2_{\mathscr{L}^2(\mathbb{R})} \lesssim ||f||^2_{\mathscr{L}^2(\mathbb{R})}$. Therefore, $||u_k||_{\mathscr{H}^1(\mathbb{R})} \lesssim ||f||_{\mathscr{L}^2(\mathbb{R})}$. Equation (2.8) follows from **proposition 1**. □

**Remark 6.** *The regularity results from this session also say that the spectrum of the operators $\mathscr{L}_{\pm k^*}$ in the bifurcation hypotheses (H1) can also be considered in $\mathscr{H}^2$ sense; therefore, the spectral splitting stated there for $\mathscr{L}^2$ still holds for $\mathscr{H}^2$.*

**Remark 7.** *Although the domain of $\mathscr{L}_k$ is $\mathscr{H}^1(\mathbb{R}; \mathbb{R}^n)$, the regularity results from the present session allow us to work on $\mathscr{H}^2(\mathbb{R}; \mathbb{R}^n)$, a Sobolev space in which these operators are bounded. Throughout this article we will work on the later rather than the first, since we want to apply some techniques from calculus in Banach spaces.*

## 3. Some auxiliary spaces

A nice consequence of lemma 1 is the characterization of any $u \in \mathscr{H}^2(\Omega; \mathbb{R}^n)$ in terms of the family $(\Pi_k[u])_{(k \in \mathbb{Z})}$. Based on this result, we will make a few definitions that will help us to better understand our problem.

**Definition 2.** *Let* $w := (\ldots, w_{-2}, w_{-1}, w_0, w_1, w_2, \ldots) \in \bigoplus_{k \in \mathbb{Z}} \mathscr{H}^2(\mathbb{R}; \mathbb{R}^n)$. *Define the norm*

$$||w||_M^2 = \sum_{k \in \mathbb{Z}} \left\{ (k^4 \vee 1)||w_k||_{\mathscr{L}^2(\mathbb{R})}^2 + (k^2 \vee 1)||\partial_x w_k||_{\mathscr{L}^2(\mathbb{R})}^2 + ||\partial_x^2 w_k||_{\mathscr{L}^2(\mathbb{R})}^2 \right\},$$

*and* $M := \{ w \in \bigoplus_{k \in \mathbb{Z}} \mathscr{H}^2(\mathbb{R}, \mathbb{R}^n) : ||w||_M < \infty \}$. *Then $M$ is a Hilbert space[6] with inner product*

$$\langle u, v \rangle_M := \sum_{k \in \mathbb{Z}} \left\{ (k^4 \vee 1)\langle u_k, v_k \rangle_{\mathscr{L}^2(\mathbb{R})} + (k^2 \vee 1)\langle \partial_x u_k, \partial_x v_k \rangle_{\mathscr{L}^2(\mathbb{R})} + \langle \partial_x^2 u_k, \partial_x^2 v_k \rangle_{\mathscr{L}^2(\mathbb{R})} \right\}.$$

**Theorem 3.1.** *Let* $S[u] := (\Pi_k[u])_{(k \in \mathbb{Z})}$, *where* $u \in \mathscr{H}^2(\Omega; \mathbb{R}^n)$. *Then,* $S : \mathscr{H}^2(\Omega; \mathbb{R}^n) \to M$. *In addition to that, $S$ is a homeomorphism between these spaces.*

*Proof.* We make use of lemma 2: as $||u||_{\mathscr{H}^1(\Omega)}^2 = \sum_{|\alpha| \leq 2} ||\partial^\alpha u||_{\mathscr{L}^2(\mathbb{R})}^2$,

$$||u||_{\mathscr{H}^1(\Omega)}^2 = \sum_{k \in \mathbb{Z}} \left\{ (k^4 + k^2 + 1)||w_k||_{\mathscr{L}^2(\mathbb{R})}^2 + (k^2 + 1)||\partial_x w_k||_{\mathscr{L}^2(\mathbb{R})}^2 + ||\partial_x^2 w_k||_{\mathscr{L}^2(\mathbb{R})}^2 \right\}$$

from which we obtain that $||S[u]||_M^2 \leq ||u||_{\mathscr{H}^2(\Omega)}^2 \leq 2||S[u]||_M^2$, which, allied to its linearity, gives the injectivity of the mapping $S(\cdot)$. It remains to be shown that it is also a surjection, so let $(w_k) \in M$. We want to find $w \in \mathscr{H}^2(\Omega; \mathbb{R}^n)$ such that $S(w) = (w_k)_{(k \in \mathbb{Z})}$; we claim that $w := \sum_{k \in \mathbb{Z}} w_k(x) e^{iky}$ is this function,[7] whose existence will be provided by the Riesz representation theorem ( which can be applied, since $M$ is a Hilbert space). Take, $\phi \in \mathscr{C}_c^\infty(\mathbb{R})$, $|\alpha| \leq 2$ and define[8]

$$T\left( \frac{\partial^\alpha}{\partial x^\alpha} \bar{\phi} \right) := \lim_{N \to \infty} \sum_{|k| \leq N} \int_{\mathbb{R}} \int_{\mathbb{T}} w_k(x) e^{iky} \cdot \frac{\partial^\alpha}{\partial x^\alpha} \bar{\phi}(x, y) dx dy.$$

In order to apply the Riesz representation theorem we need to show that T is a bounded operator in $\mathscr{H}^2(\mathbb{R}, \mathbb{R}^n)$. In details: for any fixed $N \in \mathbb{N}$,

---

[6]This is not hard to prove; the Hilbert Space property follows from the paralelogram law (see [TL86]), enhited from the (Hilbert) space $\mathscr{L}^2$, which is easily verified .

[7] This limit is taken in the $\mathscr{L}^2$ sense.

[8] We know that this sum is bounded thanks to lemma 2.

$$\sum_{|k|\leq N}\int_{\mathbb{R}}\int_{\mathbb{T}} w_k(x)e^{iky}\cdot\frac{\partial^\alpha}{\partial x^\alpha}\bar\phi(x,y)dxdy = \sum_{|k|\leq N}\int_{\mathbb{R}} w_k(x)\cdot \overbrace{(\partial^\alpha\bar\phi)_k(x)}^{\frac{1}{2\pi}\int_{\mathbb{T}}e^{iky}\frac{\partial^\alpha}{\partial x^\alpha}\bar\phi(x,y)dy}dx \overbrace{\leq}^{\text{Holder's ineq.}}$$

$$\leq \sum_{|k|\leq N}||w_k||_{\mathscr{L}^2(\mathbb{R})}||\partial^\alpha\bar\phi_k||_{\mathscr{L}^2(\mathbb{R})} \overbrace{\leq}^{\text{Cauchy's Schwartz}} \left(\sum_{|k|\leq N}||w_k||^2_{\mathscr{L}^2(\mathbb{R})}\right)^{\frac{1}{2}}\left(\sum_{|k|\leq N}||\partial^\alpha\bar\phi_k||^2_{\mathscr{L}^2(\mathbb{R})}\right)^{\frac{1}{2}} \leq$$

$$\leq ||w||_M \left(\sum_{k\in\mathbb{Z}}||\partial^\alpha\bar\phi_k||^2_{\mathscr{L}^2(\mathbb{R})}\right)^{\frac{1}{2}}.$$

As $\sum_{k\in\mathbb{Z}}||\phi_k||^2_{\mathscr{L}^2(\mathbb{R})} = ||\phi||^2_{\mathscr{L}^2(\Omega)}$, we get $\left|T\left(\frac{\partial^\alpha}{\partial x^\alpha}\phi\right)\right| \leq C||\phi||^2_{\mathscr{H}^2(\Omega)}$. Extending the operator using Hahn-Banach (or by continuity, since the $\mathscr{C}^\infty_c(\Omega)$ is dense in $\mathscr{H}^2(\Omega;\mathbb{R}^n)$ ), we prove that $T(\cdot)$ is a continuous linear functional in $M$; therefore, there exists a $w \in \mathscr{H}^2(\Omega;\mathbb{R}^n)$ such that

$$T(\phi) = \int_{\mathbb{R}}\int_{\mathbb{T}} w(x,y)\cdot\bar\phi(x,y)dxdy.$$

It remains to prove that this $w$ has the properties we are looking for; we claim that $\Pi_k(w) = w_k$, for all $k \in \mathbb{Z}$. To prove that, choose $\phi \in \mathscr{C}^\infty_c(\mathbb{R})$. We need to show that

$$T(\phi(x)e^{-iky}) = \int_{\mathbb{R}} w_k(x)\cdot\bar\phi(x)dx.$$

Fix $k \in \mathbb{Z}$. Notice that for all $N \geq |k|$,

$$\sum_{|j|\leq N}\int_{\mathbb{R}}\int_{\mathbb{T}}(w_j(x)e^{ijy})\cdot\bar\phi(x)e^{-iky}dydx = \int_{\mathbb{R}}\bar\phi(x)\cdot\left\{\sum_{|j|\leq N} w_j(x)\underbrace{\int_{\mathbb{T}}(e^{ijy})e^{-iky}dy}_{2\pi\delta_{jk}}\right\}dx =$$

$$= 2\pi\int_{\mathbb{R}}\bar\phi(x)w_k(x)dx.$$

As $N \to \infty$, we get (3). Now, since this is valid for an arbitrary $\phi$, we get the result, i.e., $\frac{1}{2\pi}\int_{\mathbb{T}} w(x,y)e^{-iky}dy = w_k(x)$ a.e., and we are done. $\square$

It worths noticing the following on the $S:\mathscr{H}^2(\Omega;\mathbb{R}^n)\to M$, proved in **appendix C**:

**Corollary 2.** *$S$ is a bounded, invertible mapping from $\mathscr{H}^2(\Omega;\mathbb{R}^n)$ to $M$. Furthermore, $S$ is a $\mathscr{C}^\infty$-diffeomorphism between these spaces.*

**Definition 3.** *Let $M_{(k_1,\ldots,k_j)} \subset M$, $w \in M_{(k_1,\ldots,k_j)}$ denote all the $w \in M$ such that $w_i(x) = 0$ whenever $i \in (k_1,\ldots,k_j)$. Denote by $||\cdot||_{M_{(k_1,\ldots,k_j)}} = ||\cdot||_M$ the norm of $M$ restricted to $M_{(k_1,\ldots,k_j)}$. For every finite set $F \subset \mathbb{Z}$, $M_{(F)}$ is a closed subspace of $M$. Following the theorem 3.1, define $S_{(k_1,\ldots,k_j)}[u] := (u_k)_{k\in\mathbb{Z}\setminus\{k_1,\ldots,k_j\}}$, for $u \in \mathscr{H}^2(\Omega;\mathbb{R}^n)$.*

Whenever $k \in \mathbb{Z}^*$ we can write (1.11) as

$$u_k = (\mathscr{L}^{\bar d}_k)^{-1}\circ\mathscr{N}_k[u] + i(\bar d - d)k(\mathscr{L}^{\bar d}_k)^{-1}[u_k]. \tag{3.2}$$

On setting $\bar{d} = d = 0$ we obtain the same system as we have in the O(2)-symmetric case. In fact, this is the main idea of the proof of both cases: we will treat the y-translational symmetric case as a perturbation of the O(2)-symmetric one; it will be a consequence of the abstract lemma on perturbations of a contraction mapping by a "small" bounded operator (see proposition (7)). The norm of this perturbation will be controled by the parameter $|\bar{d} - d|$, which is physically expected to be small. We rewrite the systems (3.2) for $k \in \mathbb{Z}^*$ as

$$S_{(0,\pm k^*)}[u] = \left((\mathscr{L}_k^{\bar{d}})^{-1} \circ \mathscr{N}_k(u_0 \oplus u_{\pm k^*} \oplus S_{(0,\pm k^*)}[u])\right)_{k \in \mathbb{Z}^*}.$$

Namely, we are looking for $w \in M_{(0,\pm k^*)}$ such that

$$w = \left((\mathscr{L}_k^{\bar{d}})^{-1} \circ \mathscr{N}_k(u_0 \oplus u_{\pm k^*} \oplus w)\right)_{k \in \mathbb{Z}^*}. \tag{3.3}$$

The solution u for the problem (1.11) will be given by setting $u = S^{-1}[u_0 \oplus u_{\pm k^*} \oplus w]$. Now that we reformulated the problem, we are ready to attack it. By hypothesis, we have that $\mathscr{L}_k$ are invertible for all $k \in \mathbb{Z}^*$. Thus, since the mapping $S[\cdot]$ provides a $\mathscr{C}^\infty$ diffeomorphism between $\mathscr{H}^2(\Omega; \mathbb{R}^n)$ and $M$ one could expect that, if we fix $u_0$ and $u_{k^*}$ in $\mathscr{H}^2(\mathbb{R}, \mathbb{R}^n)$ we would be able to solve the system (1.11) uniquely. In fact, it turns out to be the case, as we will in the next session by a contraction mapping argument/implicit function theorem.

3.1. **The contraction mapping theorem in the space** $(M_{(0,\pm k^*)}, ||\cdot||_{M_{(0,\pm k^*)}})$. From previous considerations, we must prove the following

**Theorem 3.4.** *[O(2) symmetric case] Considering the system (3.3), we define the mapping*

$$\begin{aligned} T: \quad M_{(0,\pm k^*)} &\to \bigoplus_{k \in \mathbb{Z}} \mathscr{H}^2(\mathbb{R}, \mathbb{R}^n) \\ w &\mapsto \left((\mathscr{L}_k^{\bar{d}})^{-1} \circ \mathscr{N}_k[u_0 \oplus u_{\pm k^*} \oplus w]\right)_{k \in \mathbb{Z}^*} \end{aligned}$$

*Then there exists a unique $u \in \mathscr{H}^2(\Omega; \mathbb{R}^n)$ to (1.11) such that*

$$\Pi_0(u) = u_0 \quad , \quad \Pi_{-k^*}(u) = u_{-k^*} \quad , \quad \Pi_{k^*}(w) = u_{k^*} \quad \text{and} \quad S_{(0,\pm k^*)}(u) = w.$$

*Proof.* As we are in the O(2) setting, we omit the upper index "$\bar{d}$" for now. The proof relies on the contraction mapping theorem, which will be performed first in the space $M_{(0,\pm k^*)}$. Initially, we estimate $||\mathscr{N}_k(w) - \mathscr{N}_k(v)||_{\mathscr{L}^2(\Omega)}$:

$$\begin{aligned} \mathscr{N}_k(w) - \mathscr{N}_k(v) &= [\mathscr{R}^{(1)}(v,\bar{u}) - \mathscr{R}^{(1)}(w,\bar{u})]_x + [\mathscr{R}^{(2)}(v,\bar{u}) - \mathscr{R}^{(2)}(w,\bar{u})]_y = \\ &= [df^{(1)}(v+\bar{u}) - df^{(1)}(\bar{u})]v_x - [df^{(1)}(w+\bar{u}) - df^{(1)}(\bar{u})]w_x + \\ &\quad + [df^{(1)}(v+\bar{u}) - d^2f^{(1)}(\bar{u})(v,\cdot) - df^{(1)}(w+\bar{u}) + d^2f^{(1)}(\bar{u})(w,\cdot)]\bar{u}_x + \\ &\quad + [df^{(2)}(v+\bar{u}) - df^{(2)}(\bar{u})]v_y - [df^{(2)}(w+\bar{u}) - df^{(2)}(\bar{u})]w_y \\ &= I + II + III \end{aligned}$$

Whence,

$$\begin{aligned} I &= [df^{(1)}(v+\bar{u}) - df^{(1)}(\bar{u})](v_x - w_x) + [df^{(1)}(v+\bar{u}) - df^{(1)}(w+\bar{u})]w_x \\ \Rightarrow |I| &\lesssim |v| \cdot |(v-w)_x| + |w_x| \cdot |(v-w)| \end{aligned}$$

Analogously,

$$\begin{aligned} III &= [df^{(2)}(v+\bar{u}) - df^{(2)}(\bar{u})](v_y - w_y) + [df^{(2)}(v+\bar{u}) - df^{(2)}(w+\bar{u})]w_x \\ \Rightarrow |III| &\lesssim |v| \cdot |(v-w)_y| + |w_y| \cdot |(v-w)| \end{aligned}$$

For the case II, we define $g(y) := df^{(1)}(y + \bar{u})\bar{u}_x$. We can write II as
$$II = (g(v) - g'(0)v) - (g(w) - g'(0)w) \quad \Rightarrow \quad |II| \lesssim |(v-w)|^2$$
As $\mathcal{N}_k(0) = 0$, we get

$$|T(w)|_{M_{(0,\pm k^*)}} \leq C \cdot \left(||S^{-1}[u_0 \oplus u_{\pm k^*} \oplus w]||^2_{\mathcal{H}^2(\Omega)}\right) \overset{\text{theorem 3.1}}{\leq}$$
$$C \cdot (||u_0||^2_{\mathcal{H}^2(\mathbb{R})} + ||u_{k^*}||^2_{\mathcal{H}^2(\mathbb{R})} + ||u_{-k^*}||^2_{\mathcal{H}^2(\mathbb{R})} + ||w||^2_{M_{(0,\pm k^*)}})$$

Taking $||w||_{M_{(0,\pm k^*)}} \leq R$ and $||u_0||^2_{\mathcal{H}^2(\mathbb{R})} + ||u_{k^*}||^2_{\mathcal{H}^2(\mathbb{R})} + ||u_{-k^*}||^2_{\mathcal{H}^2(\mathbb{R})} \leq \delta$, we have
$$||T[w]||_{M_{(0,\pm k^*)}} \leq C \cdot \delta + (C \cdot R) \cdot R.$$

Choosing R small in such a way that $C \cdot R \leq \theta < 1$ and $\delta > 0$ satisfying $\dfrac{C \cdot \delta}{1 - \theta} \leq R$ we get that $B := \{w \in M_{(0,\pm k^*)} : ||w||_{M_{(0,\pm k^*)} \leq R}\}$ has the property that $T(B) \subset B$. We claim that T is a contraction on B for R sufficiently small: upon estimating $||T(w) - T(v)||_{M_{(0,\pm k^*)}}$, for w, v in B and the other parameters fixed, notice that the linearity of the $\mathscr{L}_k's$ makes possible to write T(w) - T(v) as

$$T(w) - T(v) = \left(\mathscr{L}_k^{-1} \circ (\mathcal{N}_k[u_0 \oplus u_{\pm k^*} \oplus (w)] - \mathcal{N}_k[u_0 \oplus u_{\pm k^*} \oplus (v)])\right)_{k \in \mathbb{Z}^*}$$

Therefore, using the estimates on the section 2.4 and choosing R sufficiently small we obtain the result, since

$$||\mathcal{N}_k(w) - \mathcal{N}_k(v)||_{\mathscr{L}^2(\Omega)} \lesssim (||v||_{\mathscr{L}^2(\Omega)} + ||w_y||_{\mathscr{L}^2(\Omega)} + ||v - w||_{\mathscr{L}^2(\Omega)}) \cdot ||v - w||_{\mathcal{H}^2(\Omega)} \lesssim$$
$$\lesssim 2 \cdot (||v||_{\mathcal{H}^2(\Omega)} + ||w||_{\mathcal{H}^2(\Omega)}) \cdot ||v - w||_{\mathcal{H}^2(\Omega)} \lesssim (4 \cdot R)||v - w||_{\mathcal{H}^2(\Omega)}$$

Finding w as a fixed point in $M_{(0,\pm k^*)}$ we can "reassemble" u by using the mapping $S^{-1}(\cdot)$. □

As claimed before, we can bootstrap the result from the O(2)-symmetric case to the y-translational case:

**Proposition 2.** *The mapping* $(k\mathscr{L}_k^{-1}[\cdot])_{k \in \mathbb{Z}^*} : M_{(0,\pm k^*)} \to M_{(0,\pm k^*)}$ *is bounded.*

*Proof.* Denote $||\cdot||_{\mathscr{L}^2(\mathbb{R})}$ by $||\cdot||$ unless indicated; furthermore, denote $(\mathscr{L}_k^{\bar{d}})^{-1}$ by $\mathscr{L}_k^{-1}$. By definition of the norm $||\cdot||_M$,

$$||(k\mathscr{L}_k^{-1}[w])||^2_{M_{0,\pm k^*)}} = \sum_{k \in \mathbb{Z}^*} \left\{(k^4 \vee 1)||k\mathscr{L}_k^{-1}[w]||^2 + (k^2 \vee 1)||k\partial_x\mathscr{L}_k^{-1}[w]||^2 + ||k\partial_x^2\mathscr{L}_k^{-1}[w]||^2\right\}.$$

Notice that, for any $k \in \mathbb{Z}^*$

$$(k^4 \vee 1)||(k\mathscr{L}_k^{-1}[w_k]||^2 = (k^6 \vee 1)||\mathscr{L}_k^{-1}[w_k]||^2 \leq \quad (k^6 \vee 1)\frac{C}{k^4}||w_k||^2 \leq C(k^2 \vee 1)||w_k||^2,$$
$$(k^2 \vee 1)||k\partial_x\mathscr{L}_k^{-1}[w_k]||^2 \leq \quad C(k^2 \vee 1)||w_k||^2,$$
$$||k\partial_x^2\mathscr{L}_k^{-1}[w_k]||^2 \leq \quad C(k^2 \vee 1)||w_k||^2.$$

Using Parseval's lemma and **lemma 1**, $\displaystyle\sum_{k \in \mathbb{Z}^*}(k^2 \vee 1)||w_k||^2 \leq ||w_k||^2_{M_{(0,\pm k^*)}} \leq ||w||^2_{\mathcal{H}^2(\Omega)}$. Hence,

$||(k\mathscr{L}_k^{-1}[w])||_{M_{0,\pm k^*)}} \leq C||w||_{M_{(0,\pm k^*)}}.$ □

It is not hard to see that theorem 3.1 holds in the case that $\bar{d} \neq 0$, substituting $\mathscr{L}_k$ by $\mathscr{L}_k^{\bar{d}}$, *ipsis literis*. By the discussion before the last proposition, we derive the following result:

**Corollary 3.** *[Translation invariant case] Considering the system* (1.11), *and $T$ defined as in 3.1. Define the mapping*

$$T^d : \quad M_{0,\pm k^*} \quad \to \quad \bigoplus_{k \in \mathbb{Z}} \mathscr{H}^2(\mathbb{R}, \mathbb{R}^n)$$
$$w \quad \mapsto \quad T[w] + i(\bar{d} - d)((k\mathscr{L}_k^{\bar{d}})^{-1}[u_k])_{k \in \mathbb{Z}^*}$$

*Then there exists a unique $u \in \mathscr{H}^2(\Omega; \mathbb{R}^n)$ solving* (1.11) *such that*

$$\Pi_0(u) = u_0 \quad, \quad \Pi_{-k^*}(u) = u_{-k^*} \quad, \quad \Pi_{k^*}(w) = u_{k^*} \quad \text{and} \quad S_{(0,\pm k^*)}(u) = w.$$

On both cases, theorem 3.1 and corollary 3, the differentiability of the mapping $(\epsilon, u_{-k^*}, u_0, u_{k^*}) \mapsto w(\epsilon, u_{-k^*}, u_0, u_{k^*}) \in M_{(0,\pm k^*)}$, is a consequence of the parametric contraction mapping theorem (**appendix A**) and the differentiability of the mapping as studied in **appendix C**.

**Corollary 4.** *Let $w(\epsilon, u_0, u_{\pm k^*}) : \mathcal{I} \times \times_{k=1,2,3} \mathscr{H}^2(\mathbb{R}; \mathbb{R}^n) \to M_{(0,\pm k^*)}$ be the solution provided by* **theorem 3.4** *(resp corollary 3). Then*

$$w \in \mathscr{C}^1(\mathcal{I} \times \times_{k=1,2,3} \mathscr{H}^2(\mathbb{R}, \mathbb{R}^n); M_{(0,\pm k^*)}).$$

*Proof.* This is a consequence of the implicit function theorem; regularity comes from theorem C.1. □

**Remark 8.** *In virtue of the precedent results, we will denote any function $f(u)$ by $f(u_0, u_{\pm k^*})$. For instance, $\mathcal{N}_{k^*}(u)$ will be denoted by $\mathcal{N}_{k^*}(u_0, u_{\pm k^*})$.*

3.2. **Further reductions: the $k = 0$ mode.** Our first result towards the bifurcation is concerned with solving the system $\mathscr{L}_O[u_0] = u_0'' - [Df[\bar{u}]u_0]' = -\Pi_0[\partial_x(\mathscr{R}^{(1)}(u_0, \bar{u}))]$ having the other terms as parameters. So far, we have been left with the following systems

$$\mathscr{L}_0[u_0] = \mathcal{N}_0[u_0, u_{\pm k^*}] \tag{3.5a}$$
$$\mathscr{L}_{\pm k^*}[u_{\pm k^*}] = \mathcal{N}_{\pm k^*}[u_0, u_{\pm k^*}] \tag{3.5b}$$

Upon solving (3.5a) we have the equation $\mathscr{L}_O[u_0] = u_0'' - [Df(\bar{u})u_0]' = -\Pi_0[\partial_x(\mathscr{R}^{(1)}(u_0, \bar{u}))]$. Integrating from $-\infty$ to x we get[9]

$$\tilde{\mathscr{L}}_0[u_0] ; = u_0' - df(\bar{u})u_0 = \Pi_0[\mathscr{R}^{(1)}[u_0, u_{\pm k^*}]] =: \tilde{\mathcal{N}}_0[u_0, u_{\pm k^*}]. \tag{3.6}$$

However, it is well known (see [Sat76, section 2] ) that $\bar{u}_x(\cdot)$ is in the kernel of the linearized operator $\mathscr{L}^d[\cdot]$; consequently, we cannot expect $\tilde{\mathscr{L}}_0[\cdot]$ to be invertible, because $\mathscr{L}_0[\cdot]$ has a nontrivial kernel:

**Observation 4.** $\Pi_0(\bar{u}_x) = \bar{u}_x \in Ker(\mathscr{L}_0)$. *Furthermore,* $\Pi_k(\bar{u}_x) = 0$, $\forall k \neq 0$.

*Proof.* by assumption, $\bar{u}(x,y) = \bar{u}(x)$. The results in **appendix B** show that $0 = \Pi_0 \circ \mathscr{L}[\bar{u}_x] = \mathscr{L}_0[\Pi_0(\bar{u}_x)]$. Hence, $\Pi_k(\bar{u}_x) = 0$, $\forall k \neq 0$, since $\dfrac{1}{2\pi}\int_{\mathbb{T}} \bar{u}_x e^{iky} dy = \bar{u}_x \dfrac{1}{2\pi}\int_{\mathbb{T}} e^{iky} dy = \bar{u}_x \delta_{k0}$. It shows that $kernel(\mathscr{L}_0) \subset kernel(\tilde{\mathscr{L}}_0)$. □

By the regularity result proved in 3, if $v \in kernel(\tilde{\mathscr{L}}_0)$ then $v \in \mathscr{H}^2(\mathbb{R}; \mathbb{R}^n)$; we can differentiate $\tilde{\mathscr{L}}_0$, obtaining then that $v \in \tilde{\mathscr{L}}_0$, proving the set equality $kernel(\mathscr{L}_0) = kernel(\tilde{\mathscr{L}}_0)$. The operator $\tilde{\mathscr{L}}_0$ has a nontrivial kernel as well, so it cannot be invertible. We can be overcome this difficulty with the aid of an artificial constraint, allowing us to construct a right inverse $\tilde{\mathscr{L}}_0^{\dagger}[\cdot]$ (a similar treatment is developed in [MZ09]); this is the content of section 3.3, in which the following result is proved:

---
[9]As $u \in \mathscr{H}^2$, $u$, $\partial_x u \to 0$ as $|x| \to \infty$, then $u_0, u_{\pm k^*} \to 0$ as well, and then no integration constants appear in (3.6).

**Claim 1.** *Given $f \in \mathscr{L}^2(\mathbb{R}; \mathbb{R}^n)$ there exists a unique $u_0 \in \mathscr{H}^1(\mathbb{R}; \mathbb{R}^n)$ such that*

$$\tilde{\mathscr{L}}_0[u_0] = f \quad \text{and} \quad \langle u_0(0), \bar{u}_x(0) \rangle = 0 \tag{3.7}$$

**Remark 9.** *Notice, by Sobolev's embedding lemma (see [Yos95, section VI.7] ), $x \mapsto u_0(\epsilon, x) \in \mathscr{C}(\mathbb{R}; \mathbb{R}^n)$, so the evaluation of $u_0(\cdot)$ at $x = 0$ in (3.7) is legitimate. Furthermore, bootstrapping regularity with result (3) and using **corollary 6** we can obtain the following bound on the implicit functions, since $||T[u_0, u_{\pm k*}]||_{\mathscr{H}^2(\mathbb{R})} \leq C(||u_{-k*}||^2_{\mathscr{H}^2(\mathbb{R})} + ||u_0||^2_{\mathscr{H}^2(\mathbb{R})} + ||u_{k*}||^2_{\mathscr{H}^2(\mathbb{R})})$:*

$$||u_0(u_{\pm k*})||_{\mathscr{H}^2(\mathbb{R})} \leq C(||u_{-k*}||^2_{\mathscr{H}^2(\mathbb{R})} + ||u_{k*}||^2_{\mathscr{H}^2(\mathbb{R})})$$

The most important result of this section is a consequence of claim 1:

**Proposition 3.** *For every $u_{\pm k*} \in \mathscr{H}^2(\mathbb{R}; \mathbb{R}^n)$, $||u_{\pm k*}||_{\mathscr{H}^2(\mathbb{R})}$ sufficiently small, and $\epsilon \in \mathcal{I}$, there exists a unique solution $u_0(\epsilon, u_{\pm k*}) \in \mathscr{H}^2(\mathbb{R}; \mathbb{R}^n)$ to*

$$\tilde{\mathscr{L}}_0[u_0]; = u_0' - df(\bar{u})u_0 = \Pi_0[\mathscr{R}^{(1)}[u_0, u_{\pm k*}]] =: \tilde{\mathcal{N}}_0[u_0, u_{\pm k*}].$$

*satisfying $\langle u_0(0), \bar{u}_x(0) \rangle = 0$. Furthermore, $u_0 \in \mathscr{C}^1(\mathcal{I} \times \times_{k=1,2} \mathscr{H}^2(\mathbb{R}; \mathbb{R}^n); \mathscr{H}^2(\mathbb{R}; \mathbb{R}^n))$.*

Before proving this proposition, recall from **section 5** that $\mathcal{N}(w, \bar{u}) = -[\mathscr{R}^{(1)}(w, \bar{u})]_x - [\mathscr{R}^{(2)}(w, \bar{u})]_y$. We want to show that, when k = 0,

$$\Pi_0(\mathcal{N}) = -\Pi_0[\mathscr{R}^{(1)}(w, \bar{u})_x] = \partial_x[\Pi_0(\mathscr{R}^{(1)}(w, \bar{u})].$$

For the first equality it suffices to show that $\Pi_0[\mathscr{R}^{(2)}(w, \bar{u})_y] = 0$. But from the definition of $\Pi_0$,

$$\Pi_0[(\mathscr{R}^{(2)}(w, \bar{u}))_y] = \frac{1}{2\pi} \int_{\mathbb{T}} (\mathscr{R}^{(2)}(w, \bar{u}))_y dy \underbrace{=}_{2\pi periodic} 0.$$

For the second inequality, we can apply the claim on **lemma 1** as long as we prove that $\mathscr{R}^{(1)} \in \mathscr{H}^1(\mathbb{R} \times \mathbb{T})$.

**Claim 2.** $\mathscr{R}^{(1)}(w, \bar{u}) = [f^{(1)}(w + \bar{u}) - f^{(1)}(\bar{u}) - df^{(1)}(\bar{u})w]$ *is in $\mathscr{H}^1(\Omega; \mathbb{R}^n)$, whenever $u, w \in \mathscr{H}^2(\Omega; \mathbb{R}^n)$. Furthermore, $||\mathscr{R}^{(1)}(w, \bar{u})||_{\mathscr{H}^1(\Omega)} \lesssim ||w||^2_{\mathscr{H}^2(\Omega)}$.*

*Proof.* This is a consequence of the density of smooth, compactly supported functions in $\mathscr{H}^2(\Omega; \mathbb{R}^n)$; the estimates to show that $\mathscr{R}^{(1)}(w, \bar{u})$ and its weak derivatives are in $\mathscr{L}^2(\Omega; \mathbb{R}^n)$ follow from the differentiability of $f^{(1)}$ and $f^{(2)}$ and the embedding theorem 2.1. The proof of $||\mathscr{R}^{(1)}(w, \bar{u})||_{\mathscr{H}^1(\Omega)} \lesssim ||w||^2_{\mathscr{H}^2(\Omega)}$ is essentially the same as those in theorem 3.4. □

As we will show later (**theorem 3**), $||\tilde{\mathscr{L}}_0^\dagger[f]||_{\mathscr{H}^2(\mathbb{R})} \lesssim ||f||_{\mathscr{H}^1(\mathbb{R})}$. Formally, the solution to (3.6) is given by $u_0 = \tilde{\mathscr{L}}_0^\dagger \circ \tilde{\mathcal{N}}_0[u_0, u_{\pm k*}]$. Define the mapping $\mathscr{T}[w] := \tilde{\mathscr{L}}_0^\dagger \circ \tilde{\mathcal{N}}_0[w, u_{\pm k*}]$, where $u_{\pm k*} \in \mathscr{H}^2(\mathbb{R}; \mathbb{R}^n)$ are fixed. We need to show that for every $u_{\pm k*}$ sufficiently small in $||\cdot||_{\mathscr{H}^2(\mathbb{R})}$ norm we have a fixed point $w = \mathscr{T}[w]$ in $\mathscr{H}^2(\mathbb{R}; \mathbb{R}^n)$.

*Proof of the proposition 3.* Define $f[u_0, u_{\pm k*}] = \tilde{\mathscr{L}}_0[u_0] - \tilde{\mathcal{N}}_0[u_0, u_{\pm k*}]$, where $\frac{\partial}{\partial u_0}\tilde{\mathscr{L}}_0\Big|_{u_0=w} = \tilde{\mathscr{L}}_0[w]$, since $\tilde{\mathscr{L}}_0$ is linear. The openness of the invertible operators and the invertibility of the later operator will imply that $\frac{\partial}{\partial u_0}f(u_0,)$ is an isomorphism on $\mathscr{H}^2(\mathbb{R}; \mathbb{R}^n) \mapsto \mathscr{L}^2(\mathbb{R}; \mathbb{R}^n)$ if we show that $||\frac{\partial}{\partial u_0}\tilde{\mathcal{N}}_0[u_0, u_{\pm k*}]||_{\mathscr{H}^2}$ is sufficiently small; the argument closes with an application of the implicit function theorem in Banach spaces (see [CH82, Cra77] ). In order to derive our result, notice that $||w(u_0, u_{\pm k*})||_{M_{0, \pm k*}} = \mathcal{O}(||u_0||^2_{\mathscr{H}^2(\mathbb{R})}, ||u_{\pm k*}||^2_{\mathscr{H}^2(\mathbb{R})})$, by remark 8. We know from appendix

C that $(\epsilon, u) \mapsto \tilde{\mathcal{N}}_0[\epsilon, u]$ is a $\mathscr{C}^1(\mathcal{I} \times \mathscr{H}^2(\mathbb{R}; \mathbb{R}^n); M_{0, \pm k^*})$ mapping. Similarly to the estimates in theorem 3.4, one can show that $||\tilde{\mathcal{N}}_0[u_0 + h, \tilde{u}] - \tilde{\mathcal{N}}_0[u_0, \tilde{u}]||_{\mathscr{H}^1(\mathbb{R})} \lesssim ||h||^2_{\mathscr{H}^2(\mathbb{R})}$. As $\tilde{\mathcal{N}}_0[0, 0] = 0$ and $(\epsilon, u) \mapsto \tilde{\mathcal{N}}_0[\epsilon, u] \in \mathscr{C}^1(\mathcal{I} \times \mathscr{H}^2; \mathscr{H}^1)$ we combining the above results in order to get that the Frèchét derivative $\partial_{u_0} \tilde{\mathcal{N}}_0 \Big|_{(u_0, \tilde{u}) = (0,0)} = 0$; continuity with respect to parameters provides $||\partial_{u_0} \tilde{\mathcal{N}}_0[u_0, \tilde{u}]||$ small for $(u_0, \tilde{u})$ sufficiently small. It follows then that $\dfrac{\partial}{\partial u_0} f(u_0)$ is an isomorphism. Regularity is a consequence of the implicit function theorem (see **appendix C**). □

**Remark 10.** *In virtue of the precedent results, we will denote any function $f(u_0, u_{\pm k^*})$ by $f(u_{\pm k^*})$ only. For instance, $\mathcal{N}_{k^*}(u_0, u_{\pm k^*})$ will be denoted by $\mathcal{N}_{k^*}(u_{\pm k^*})$.*

3.3. **Proof of claim 1: the existence of $\tilde{\mathscr{L}}_0^\dagger$.** Remark 2 allow us to apply the gap/conjugation lemma (as explained in [MZ05, chapter 2]) to the problem (3.6), obtaining then a mapping $W(\cdot, \cdot) : \mathcal{I} \times \mathbb{R} \setminus \{0\} \to \mathbb{C}^{n \times n}$ that conjugates the solutions of system $\tilde{\mathscr{L}}_0[u_0] = f$ and the solutions of the asymptotic system when $|x| \to \pm \infty$ and so that $x \mapsto W^{-1}(\epsilon, x)$ is uniformly bounded in $\mathbb{R}$ and $\lim_{x \to \infty} W(\epsilon, x) \to Id$; in details, the gap lemma asserts that $W(\epsilon, x) z(x)$ solves (3.6) if and only if z solves

$$z' - A^{\pm}(\epsilon) z = W^{-1}(\epsilon, x) f = g^{\pm} \quad , \quad x \gtrless 0,$$

In this section we make a quick review, with some adjustments, of [Hen81, theorem A.1, chapter 5]: let $\mathbb{P}^+(\epsilon)$ and $\mathbb{Q}^-(\epsilon)$ be projections on the eigenspaces with positive real part of $A^+(\epsilon)$ and negative real parts of $A^-(\epsilon)$, respectively. The case $x > 0$, considered as a invitial value problem, has a solution

$$z(x) := e^{A^+(x - x_0^+)} \mathbb{P}^+ z_0^+ + \int_{x_0}^{x} e^{A^+(x-s)} \mathbb{P}^+ g(s) ds - \int_{x}^{+\infty} e^{A^+(x-s)} (I - \mathbb{P}^+) g(s) ds \quad , \quad \text{for} \quad x \geq 0.$$

for $x_0^+ > 0$. Analogously for the case $x < 0$,

$$z(x) := e^{A^-(x - x_0^-)} (I - \mathbb{Q}^-) z_0^- - \int_{x}^{x_0} e^{A^-(x-s)} (I - \mathbb{Q}^-) g(s) ds + \int_{-\infty}^{x} e^{A^-(x-s)} \mathbb{Q}^- g(s) ds \quad , \quad \text{for} \quad x < 0.$$

where $x_0^- < 0$. Taking $x_0^+ \downarrow 0$, $x_0^- \uparrow 0$, and noticing that both integrals can be seen as convolution operators we rewrite the above equations as

$$z(x) := e^{A^{\pm}(\epsilon) x} M^{\pm} z_0^{\pm} + \int_{\mathbb{R}} H^{\pm}(\epsilon, x - s) 1^{\pm}(s) g^{\pm}(s) ds \quad , \quad \text{for} \quad x \gtrless 0.$$

where $M^+(\epsilon, x) = e^{A^+(\epsilon)(x)} \mathbb{P}^+(\epsilon)$, $M^-(\epsilon, x) = e^{A^-(\epsilon)(x)} (I - \mathbb{Q}^-(\epsilon))$, $1^+ = 1_{[0, +\infty]}$ and $1^- = 1_{[-\infty, 0]}$. We need to show that one can find $z_0^{\pm} \in \mathbb{C}^n$ so that both equations above have a solution. In fact, we can solve this problem by satisfying the transmission condition $W(\epsilon, 0^+) z(0^+) = W(\epsilon, 0^-) z(0^-)$, since we want solutions to be continuous. Therefore

$$W(\epsilon, 0^+)[z_0^+] - W(\epsilon, 0^-)[z_0^-] = -\left\{ W(\epsilon, 0^+) H^+ * 1^+ g^+ - W(\epsilon, 0^-) H^- * 1^- g^- \right\}. \tag{3.8}$$

$z_0^+$ and $z_0^-$ can be found provided that $span\{W(0^+)[\mathbb{E}_+^s] \oplus W(0^-)[\mathbb{E}_-^u]\} = \mathbb{C}^N$. We know that $\mathbb{E}_+^s \oplus \mathbb{E}_-^u = \mathbb{C}^N$ and that $dim(\mathbb{E}_+^s) + dim(\mathbb{E}_-^u) = N + 1$. The only thing that we need to prove is that $W(\epsilon, 0^+)$ and $W(\epsilon, 0^-)$ preserve "transversality" of these subspaces and that these mapping have nice differentiability properties in $\epsilon$. In fact, this result follows from $dim[W(\epsilon, 0^+)[\mathbb{E}_+^s] \cap W(\epsilon, 0^-)[\mathbb{E}_-^u]] = 1$, which is true, since $Ker(\tilde{\mathscr{L}}_0) = \{\bar{u}_x^\epsilon\}$ (observation 4 ). Recall that if u solves $\mathscr{L}_0[u_0] = f$ then $u = \tilde{u} + t\bar{u}_x$ ($\tilde{u}$ is a particular solution); if there exists a vector $l$ such that $l \cdot \bar{u}_x(0) \neq 0$ but $l \cdot [\tilde{u} + t\bar{u}_x](0) = 0$, we claim that we have

**Proposition 4.** *[Uniqueness of solutions to problem (3.7)] As constraints to problem $\mathscr{L}_0[u_0] = f$, consider $l \cdot \bar{u}_x(0) \neq 0$ and $l \cdot [\tilde{u} + t\bar{u}_x](0) = 0$. Then the solutions to this problem is unique.*

*Proof.* Assume u and v are solutions. Then $\mathscr{L}_0[u] = \mathscr{L}_0[v] = f \Rightarrow \mathscr{L}_0[u - v] = 0$, hence $u - v = t\bar{u}_x \Rightarrow$. Then, $l \cdot (u - v)(0) = 0 = tl \cdot \bar{u}_x(0) \Leftrightarrow t = 0$. It follows that $u = v$. □

We investigate now the consequences of condition $\langle u(0), \bar{u}_x(0) \rangle = 0$ on the compatibility conditions asserted by equation (3.8). Let $r = dim[W(\epsilon, 0^+)[\mathbb{E}_+^s]$ and $s = dim[W(\epsilon, 0^-)[\mathbb{E}_-^u]]$, $r + s = N + 1$ and $span\{v_1, \ldots, v_r\} = W(\epsilon, 0^+)[\mathbb{E}_+^s]$, $span\{w_1, \ldots, w_s\} = W(0^-)[\mathbb{E}_-^u]$ basis. Then

$$W(0^+)z_0^+ = [v_1 \ldots v_r]\alpha \quad , \quad \alpha \in r \times 1 \quad \text{and} \quad W(0^-)z_0^- = [w_1 \ldots w_s]\beta \quad , \quad \beta \in s \times 1$$

Using the Lax shock condition (H3), $\bar{u}_x^\epsilon$ is so that $\bar{u}_x^\epsilon(0) = [v_1 \ldots v_r](\epsilon)\alpha^*(\epsilon) = [w_1 \ldots w_s](\epsilon)\beta^*(\epsilon)$. We end up with the following matrix $M(\epsilon) \in \mathbb{C}^{(N+1)\times(N+1)}$

$$\underbrace{\begin{bmatrix} v_1 \cdots v_r & w_1 \cdots w_s \\ \alpha^{*\prime} & 0 \end{bmatrix}}_{M(\epsilon)} \begin{bmatrix} \alpha \\ \beta \end{bmatrix} = \begin{bmatrix} G \otimes g \\ 0 \end{bmatrix}$$

As all the entries are $\mathscr{C}^k(\mathcal{I}; \mathbb{C})$, $M(\epsilon, \cdot)$ is $\mathscr{C}^k$. It remains to prove that $M(\epsilon)$ is invertible. It suffices to prove that

$$M(0)\begin{bmatrix} \alpha \\ \beta \end{bmatrix} = 0.$$

Thus, $[v_1 \ldots v_r]\alpha = -[w_1 \ldots w_s]\beta = \tilde{u} \Rightarrow \tilde{u} \in W(0^+)[\mathbb{E}_+^s] \cap W(0^-)[\mathbb{E}_-^u]$. Then $\alpha = t\alpha^*$. Considering only the last line of the matrix $M(\epsilon)$, we have

$$\begin{bmatrix} \alpha^* & 0 \end{bmatrix} \begin{bmatrix} t\alpha^* \\ \beta \end{bmatrix} = t\langle \alpha^*, \alpha^* \rangle = 0.$$

As $\alpha^* \neq 0$, we have that $t = 0$, thus $\alpha = 0$. But then $[w_1 \ldots w_s]\beta = 0$. As $\{w_1 \ldots w_s\}$ are linearly independent, $\beta = 0$. Therefore, M is invertible. Therefore, we can write the solution as $u(\epsilon, x) = \tilde{H}(\epsilon) * f$. As $\partial_\epsilon \tilde{H}(\epsilon)(\cdot) \in \mathscr{L}^2(\mathbb{R}; \mathbb{R}^n)$ for all $j \in \{1, 2, \ldots k\}$ and $\tilde{\mathscr{L}}_0$ is linear, we have that $u(\epsilon, f) \in \mathscr{C}^1(\mathcal{I} \times \mathscr{L}^2(\mathbb{R}; \mathbb{R}^n); \mathscr{L}^2(\mathbb{R}; \mathbb{R}^n))$. Using convolution inequalities (namely, $||f * g||_{\mathscr{L}^\infty} \leq ||f||_{\mathscr{L}^2}||g||_{\mathscr{L}^2}$ and $||f * g||_{\mathscr{L}^2} \leq ||f||_{\mathscr{L}^1}||g||_{\mathscr{L}^2}$ we obtain

$$||z||_{\mathscr{L}^2(\mathbb{R})} \leq C(\epsilon)(||g^+||_{\mathscr{L}^2(\mathbb{R}^+)} + ||g^-||_{\mathscr{L}^2(\mathbb{R}^-)}).$$

To complete the argument, we prove a regularity lemma that - although not given in [Hen81] - is implicitly assumed:

**Claim 3.** *If the system (3.6) is solvable in $\mathscr{L}^2(\mathbb{R})$, i.e., $||u_0||_{\mathscr{L}^2(\mathbb{R})} \lesssim ||f||_{\mathscr{L}^2(\mathbb{R})}$ then $u_0 \in \mathscr{H}^1(\mathbb{R})$ and $||u_0||_{\mathscr{H}^1(\mathbb{R})} \lesssim ||f||_{\mathscr{L}^2(\mathbb{R})}$. In case $f \in \mathscr{H}^1(\mathbb{R})$, we have that $u_0 \in \mathscr{H}^2(\mathbb{R})$ and $||u_0||_{\mathscr{H}^2(\mathbb{R})} \lesssim ||f||_{\mathscr{H}^1(\mathbb{R})}$*

**Remark 11.** *By similar techniques as those presented in appendix C, this result implies that $u_0(\epsilon, f) \in \mathscr{C}^1(\mathcal{I} \times \mathscr{H}^1(\mathbb{R}; \mathbb{R}^n); \mathscr{H}^2(\mathbb{R}; \mathbb{R}^n))$.*

*Proof of claim 3.* We want to prove that, for every $\phi(\cdot) \in \mathscr{C}_0^\infty(\mathbb{R}; \mathbb{R}^n)$, we have that

$$-\int_\mathbb{R} u \cdot \phi'(x)dx = \int_\mathbb{R} (f + A[\bar{u}(x)]u(x)) \cdot \phi(x)dx \tag{3.9}$$

Take $f_n \in \mathscr{C}_0^\infty(\mathbb{R}; \mathbb{R}^n) \cap \mathscr{L}^2(\mathbb{R}; \mathbb{R}^n) \to f$ in $\mathscr{L}^2$ norm. The solutions $u_n$ to $\mathscr{L}_0[u_n] = f_n$ given through the method in **section 3.3** are then classical solutions and satisfy the integral equality

$$-\int_\mathbb{R} u_n \cdot \phi'(x)dx = \int_\mathbb{R} (f_n + A[\bar{u}(x)]u_n(x)) \cdot \phi(x)dx.$$

By linearity, we also have that $\mathscr{L}_0[u_n - u_m] = f_n - f_m$ It implies that $u_{n n \in \mathbb{N}}$ is a Cauchy sequence in $\mathscr{L}^2(\mathbb{R}, \mathbb{R}^n)$. We can find a limit point $u \in \mathscr{L}^2(\mathbb{R}; \mathbb{R}^n)$ such that $\mathscr{L}_0[u] = f$. Recall that the asymptotic limit of $\bar{u}(x)$ as $|x| \to \pm\infty$ show that $x \to ||A[\bar{u}(x)]|||$ is uniformly bounded in $\mathbb{R}$; therefore, $A[\bar{u}(\cdot)]u_n(\cdot) \to A[\bar{u}(\cdot)]u(\cdot)$ in $\mathscr{L}^2$. As strong convergence implies weak convergence, taking $n \to \infty$ on the previous equation gives that $u \in \mathscr{H}^1(\mathbb{R}; \mathbb{R}^n)$. Now, assume that $f \in \mathscr{H}^1(\mathbb{R}; \mathbb{R}^n)$. One can rewrite (3.9) as

$$-\int_\mathbb{R} u \cdot \phi''(x)dx = \int_\mathbb{R} (f + (A(\bar{u}(x)u(x)))) \cdot \phi'(x)dx$$

By the same method as before, one can show that $A(\bar{u}(\cdot))u(\cdot)$ has weak derivative in $\mathscr{L}^2(\mathbb{R}; \mathbb{R}^n)$ given by $(A[\bar{u}(\cdot)]u(\cdot))' = A'[\bar{u}(\cdot)]\bar{u}'(\cdot)u(\cdot) + A[\bar{u}(\cdot)]u'(\cdot)$. We can bootstrap our energy estimate to $\mathscr{H}^2$ using the energy estimate we just obtained for $\mathscr{H}^1$, to get

$$||u||_{\mathscr{H}^2(\mathbb{R})} \lesssim ||u||_{\mathscr{H}^1(\mathbb{R})} + ||u||_{\mathscr{H}^1(\mathbb{R})} \lesssim ||f||_{\mathscr{H}^1(\mathbb{R})}$$

$\square$

## 4. Symmetry and Lyapunov-Schmidt method

Our main concern in this part is on understanding how the infinite dimensional system (1.11) inherits the symmetries contained in the equation (1.2), since it is not clear that the mapping S respects these features of the model. It turns out that the geometry of the projections on the Fourier modes and the spectral synthesis in the y direction provide the necessary consistency between the symmetries on both cases, as we will show in this section. We split the analysis into the translation invariant case (common to both systems ) and the $O(2)$ case.

### 4.1. Robustness of symmetries through Lyapunov-Schmidt reduction: translation invariant case. 
Let $\tau_c$ be a translation operator $\tau_c[u](x, y) = u(x, y + c)$. It is clear, from the y independence of the underlying wave $\bar{u}$, that $\tau_c[w]$ is a family of solutions to (1.6) whenever $w$ solves it. Summarizing:

**Proposition 5.** *$u(\cdot, \cdot)$ solves (1.6) if and only if $\tau_c[u](\cdot, \cdot)$ solves it, $\forall c \in \mathbb{R}$. Furthermore, in case of a pair of eigenvector/eigenvalue $(u/\lambda)$ of the operator $\mathscr{L}$, $(\tau_c[u], \lambda)$ is also a eigenvector/eigenvalue pair $\forall c \in \mathbb{R}$ of the same operator.*

We focus now on understanding how the symmetry is inherited by the family of operators $\{\mathscr{L}_k\}_{k \in \mathbb{Z}}$. It turns out that the only thing that must be studied is the relation between the mappings $\Pi_k$ ( $k \in \mathbb{Z}$) and $\tau_c$ ($c \in \mathbb{R}$). In fact,

$$(\Pi_k \circ \tau_c[u])(x) = \frac{1}{2\pi}\int_\mathbb{T} \tau_c[u](x, \xi)e^{-ik\xi}d\xi = e^{ikc}\{\Pi_k[u](x)\} \quad \Rightarrow \quad \Pi_k \circ \tau_c[\cdot] = e^{ikc}\Pi_k[\cdot].$$

If we define the mapping $\tilde{\tau}_c[(u_k)]_{\{k \in F \subset \mathbb{Z}\}} \mapsto (e^{ikc}u_k)_{\{k \in F \subset \mathbb{Z}\}}$ the previous result implies that

$$\tilde{\tau}_c \circ S = S \circ \tau_c \tag{4.1}$$

in $\mathscr{H}^2(\Omega; \mathbb{R}^n)$. Thus, if $(u_k)_{k \in \mathbb{Z}}$ then $c \mapsto (e^{ikc}u_k)_{k \in \mathbb{Z}}$ is a one parameter family of solutions in the space M. An analogue of the previous proposition can be translated in terms of the $\mathscr{L}_{k k \in \mathbb{Z}}$:

**Proposition 6.** $(u_k)_{k \in \mathbb{Z}}$ *solves (1.11) if and only if* $(e^{ikc}u_k)_{k \in \mathbb{Z}}$ *solves it,* $\forall c \in \mathbb{R}$.

It is natural then to ask if this property is preserved upon the Lyapunov-Schmidt reduction we performed. In order to prove that, let $(u_0, u_{\pm k^*})$ and $w(u_0, u_{\pm k^*})$ be functions as in theorem 3.4; symmetry preserving is the same as showing that $w[(u_0, e^{i\pm k^* c}u_{\pm k^*})] = \tilde{\tau}_c w[(u_0, u_{\pm k^*}]$; this is the content of the next

**Proposition 7.** *The Lyapunov Schmidt reduction performed on section 3 respects y-translational symmetry, namely,* $w(\tau_c[u_0, u_{\pm k^*}]) = \tau_c[w(u_0, u_{\pm k^*})]$

*Proof.* Fix $c \in \mathbb{R}$. Consider the modified mapping

$$T_{\tilde{\tau}_c}: M_{(\pm k^*, 0)} \to [\mathscr{H}^2(\mathbb{R}; \mathbb{R}^n)]^3$$
$$w \mapsto \left(\mathscr{L}_k^{-1} \circ \mathscr{N}_k[\tilde{\tau}_c[u_0 \oplus u_{\pm k^*} \oplus w]]\right)_{k \in \mathbb{Z}^*}.$$

Notice that the operator $\tilde{\tau}_c$ is an isometry. So we are in the same neighborhood as in **theorem 3.4**, and then $T_{\tilde{\tau}_c}$ is a contraction, being the proof essentially the same as in 3.4. We know that $(u_0, u_{\pm k^*}, w)$ solves equation (1.11) if and only if $\tilde{\tau}_c(u_0, u_{\pm k^*}, w)$ solves it as well. From property 4.1 and uniqueness of the fixed point, we conclude that $w[\tilde{\tau}_c(u_0, u_{\pm k^*})] = \tilde{\tau}_c[(u_0, u_{\pm k^*}, w)]$. □

**4.2. Robustness of symmetries through Lyapunov-Schmidt reduction: the O(2) case.**
By assumption, if v solves (1.11) then $\Gamma[v]$ also solves it. Initially we need to understand the action of $\Gamma$ (see definition 1.5) on the space M. As before, we want to prove that rotational symmetry is inherited upon Lyapunov-Schmidt reduction.

**Proposition 8.** *The Lyapunov Schmidt reduction performed on **section 3** respects reflective symmetry, namely,* $w(\tilde{\Gamma}[u_0, u_{\pm k^*}]) = \tilde{\Gamma}[w(u_0, u_{\pm k^*})]$.

*Proof.* Define the mapping

$$T: M_{(\pm k^*, 0)} \Rightarrow [\mathscr{H}^2(\mathbb{R}; \mathbb{R}^n)]^3$$
$$w \mapsto \left(\mathscr{L}_k^{-1} \circ \mathscr{N}_k[[\Gamma[u_0 \oplus u_{k^*} \oplus u_{-k^*}] \oplus w]]\right)_{k \in \mathbb{Z}^*}$$

Notice that $\|\Gamma[u_{k^*}]\|_{\mathscr{H}^2(\mathbb{R})} = \|u_{-k^*}\|_{\mathscr{H}^2(\mathbb{R})}$, since $R \in O(n)$. By uniqueness of the contraction mapping and the fact that u solves (1.11) if and only if $\Gamma u$ solves it too, we get the result. □

**4.3. Robustness of symmetries through implicit function theorem.** As pointed out before, as k = 0 there is no difference between the two cases in which concerns to the equations; the analysis in both - translation invariant and O(2) case - is very similar then, allowing us to deal with then together. The only thing that we need to prove then is that the solution given by the right inverse in 3 preserves the O(2) symmetry. In fact, we want to show that, given a solution $u_0$ to (3.6) satisfying the phase condition $\langle u(0), \bar{u}_x(0) \rangle = 0$, then $u_0(\Gamma(u_{\pm k^*})) = \Gamma[u_0(u_{\pm k^*})]$. Define then

$$\tilde{\mathscr{L}}_0[u_0] = \tilde{\mathscr{N}}_0[u_0, \Gamma(u_{\pm k^*})].$$

$\|\Gamma(u_{\pm k^*})\|_{\mathscr{H}^2(\mathbb{R})} = \|u_{\pm k^*}\|_{\mathscr{H}^2(\mathbb{R})}$, since $R \in O(n)$. Therefore the same computations as those used in 3.3 provide existence of a solution $u_0 = u_0(\Gamma(u_{\pm k^*}))$. Now, we know by symmetry that if $(u_0, u_{\pm k^*})$ solves the system then $\Gamma(u_0, u_{\pm k^*})$ also solves it. We only need to show that it satisfies the artificial constraint $\Gamma(u_0)(0) \cdot \bar{u}_x(0) = 0$. It turns out that this is a consequence of $R \in O(n)$: R is an isometry therefore it also preserve inner products. Then, $\Gamma(u_0)(0) \cdot \bar{u}_x(0) = R(u_0)(0) \cdot \bar{u}_x(0) = (u_0)(0) \cdot R\bar{u}_x(0)$. As $R\bar{u}(x) = \bar{u}(x)$, $R\bar{u}_x(x) = \bar{u}_x(x)$. Hence,

$$\Gamma(u_0)(0) \cdot \bar{u}_x(0) = R(u_0)(0) \cdot \bar{u}_x(0) \overbrace{=}^{R\bar{u}_x(x) = \bar{u}_x(x)} R(u_0)(0) \cdot R\bar{u}_x(0) = u_0(0) \cdot \bar{u}_x(0) = 0$$

Therefore we must have $u_0(\Gamma(u_{\pm k^*})) = \Gamma[u_0(u_{\pm k^*})]$, since the solutions are unique in a neighborhood of $u_0 = 0$. It finishes the proof.

## 5. A FURTHER STEP TOWARDS THE REDUCED EQUATIONS

So far, we have reduced problem (1.11) to problem

$$\mathscr{L}_{\pm k^*}^{\bar{d}}[u_{\pm k^*}] = \mathscr{N}_{\pm k^*}[u] \pm ik^*(\bar{d} - d)u_{\pm k^*} \tag{5.1}$$

We will solve the case $k = k^*$, since the other is similarly treated. Let $\mathcal{P}_{k^*}(\epsilon, \cdot) : \mathscr{H}^2(\mathbb{R}; \mathbb{R}^n) \to \mathscr{H}^2(\mathbb{R}; \mathbb{R}^n)$ be the spectral projection onto the eigenspace associated to the $\lambda_{k^*}(\epsilon)$ of $\mathscr{L}_{k^*}^{\bar{d}}$:

$$\mathcal{P}_{k^*}(\epsilon, v) := \frac{1}{2\pi i} \int_\gamma (\mu I - \mathscr{L}_{k^*}^{\bar{d}})^{-1}[v] d\mu \tag{5.2}$$

The results in C asserts, by differentiation under the integral sign, that $(\epsilon, v) \mapsto \epsilon^{(i)} \mathcal{P}_{k^*}(\epsilon, v)$ for $i \in \{0, 1, 2\}$ exists and are bounded for $\epsilon \in \mathcal{I}$; therefore the composition with the mapping $\mathcal{P}_{k^*}(\cdot, \cdot)$ will not change regularity of the operators in equation (5.1). By property (H1) we can choose a contour $\gamma$ independent of $\epsilon$ so that the projections are well defined and have the desired properties (see remark 6). We can use $\mathcal{P}_{k^*}$ to rewrite 5.1 as

$$\mathscr{P}_{k^*}\mathscr{L}_{k^*}[u_{k^*}] = \mathscr{P}_{k^*}\mathscr{N}_{k^*}[\mathscr{P}_{k^*}\tilde{u} + (I - \mathscr{P}_{k^*})\tilde{u}] + i(\bar{d} - d)k^* \mathscr{P}_{k^*} u_{k^*} \tag{5.3a}$$
$$(I - \mathscr{P}_{k^*})\mathscr{L}_{k^*}[u_{k^*}] = (I - \mathscr{P}_{k^*})\mathscr{N}_{k^*}[\mathscr{P}_{k^*}\tilde{u} + (I - \mathscr{P}_{k^*})\tilde{u}] + i(\bar{d} - d)k^*(I - \mathscr{P}_{k^*})u_{k^*}. \tag{5.3b}$$

As $\lambda_{k^*}$ is a simple eigenvalue, $dim(Range\mathcal{P}_{k^*}) = 1$. Thus, to reduce the problem to a finite dimensional one we could try to solve **equation** (5.3b) for $(I - \mathcal{P}_{k^*})u_{k^*}$ in terms of $\mathcal{P}_{k^*}[u_{k^*}]$ and $\epsilon$, which turns out to be possible, as we show next. As we did in **section 3**, we prove the O(2)-symmetric case first and obtain the translational-symmetric case by a perturbation method.

**Proposition 9.** *[O(2) Symmetric case] Considering the equation* (5.3b) *in the O(2)-symmetric setting, we define the mapping*

$$W[v, u] := \mathscr{L}_{k^*}^{-1}(I - \mathscr{P}_{k^*})\mathscr{N}_{k^*}[u + v]$$

*where* $u \in Range(\mathcal{P}_{k^*}) = \mathcal{P}_{k^*}(\mathscr{H}^2(\mathbb{R}; \mathbb{R}^n))$ *and* $v \in Range(I - \mathcal{P}_{k^*}) = (I - \mathcal{P}_{k^*})(\mathscr{H}^2(\mathbb{R}; \mathbb{R}^n))$. *Then, the mapping* $v \mapsto W[v, u]$ *has a fixed point* $v = v(u, \epsilon)$ *for all* $||u||_{\mathscr{H}^2(\mathbb{R})}$ *small. Furthermore,* $v \in \mathscr{C}^1(\mathcal{I} \times H^2(\mathbb{R}; \mathbb{R}^n); H^2(\mathbb{R}; \mathbb{R}^n))$.

*Proof.* Initially, notice that $\mathcal{P}_{k^*}$ has a closed range. Therefore $(I - \mathcal{P}_{k^*})(\mathscr{H}^2(\mathbb{R}))$ is a closed subspace of $\mathscr{H}^2(\mathbb{R}; \mathbb{R}^n)$; in particular, it is a Banach space with the induced norm, so we can apply the the contraction mapping theorem. Since $\mathscr{P}_{k^*}(I - \mathscr{P}_{k^*}) = 0$, 5.3b can be solved and the mapping $\mathscr{L}_k^{-1}(I - \mathscr{P}_{k^*})$ is a bounded operator from $\mathscr{L}^2(\mathbb{R}; \mathbb{R}^n)$ to $(I - \mathcal{P}_{k^*})\mathscr{H}^2(\mathbb{R}; \mathbb{R}^n)$, thanks to results in appendix D. Fix $u \in \mathcal{P}_{k^*}(\mathscr{H}^2(\mathbb{R}; \mathbb{R}^n))$. Using the estimates derived in **theorem 3.4**, we get

$$||\mathscr{N}_k[v + u]] - \mathscr{N}_k[\tilde{v} + u]|| \leq \theta ||v - \tilde{v}||_{\mathscr{H}^2(\mathbb{R})}^2.$$

Following similar arguments as those in the proof of theorem 3.4 we have a contraction for $||v||_{\mathscr{H}^2(\mathbb{R})}^2$, $||\tilde{v}||_{\mathscr{H}^2(\mathbb{R})}^2$ sufficiently small. Furthermore, $W(B) \subset B$ for $||u||_{\mathscr{H}^2(\mathbb{R})}$ sufficiently small, since

$$||W(v + u)|| \leq ||\mathscr{L}_k^{-1}(I - \mathscr{P}_{k^*})||_{\mathscr{L}^2 \to \mathscr{H}^2} \cdot ||\mathscr{N}_{k^*}[\mathscr{P}_{k^*}\tilde{u} + (I - \mathscr{P}_{k^*})\tilde{u}]||_{\mathscr{L}^2(\mathbb{R})} \lesssim$$
$$\lesssim ||v + u||_{\mathscr{H}^2(\mathbb{R})}^2 \lesssim ||v||_{\mathscr{H}^2(\mathbb{R})}^2 + ||u||_{\mathscr{H}^2(\mathbb{R})}^2.$$

The differentiability follows from the differentiability of the mapping $\mathcal{N}_{k^*}$ in $\epsilon$ and $u$, as proved in the appendix C. □

**Corollary 5.** *[Translation invariant case] Considering the equation (5.3b) in the translational-symmetric setting, we define the mapping*

$$\tilde{W}[v,u] := (\mathcal{L}_{k^*}^{\bar{d}})^{-1}(I - \mathcal{P}_{k^*})\mathcal{N}_{k^*}[u+v] + i(\bar{d}-d)k^*(\mathcal{L}_{k^*}^{\bar{d}})^{-1}(I - \mathcal{P}_{k^*})(v+u)$$

*where $u \in Range(\mathcal{P}_{k^*}) = \mathcal{P}_{k^*}(\mathcal{H}^2(\mathbb{R};\mathbb{R}^n))$ and $v \in Range(I - \mathcal{P}_{k^*}) = (I - \mathcal{P}_{k^*})(\mathcal{H}^2(\mathbb{R};\mathbb{R}^n))$. Then, the mapping $v \mapsto W[v,u]$ has a fixed point $v = v(u,\epsilon)$ for all $||u||_{\mathcal{H}^2(\mathbb{R})}$ small. Furthermore, $v \in \mathscr{C}^1(\mathcal{I} \times H^2(\mathbb{R};\mathbb{R}^n); H^2(\mathbb{R};\mathbb{R}^n)))$.*

*Proof.* Just notice that the mapping $v \mapsto P(v,u) := (\bar{d}-d)k^*(v+u)$ is bounded and $P(v,u) - P(\tilde{v},u)$ depends only on the norm $||v - \tilde{v}||_{\mathcal{H}^2(\mathbb{R})}$ and the smallness of $|\bar{d}-d|$. The result follows from the **corollary 7**, which also provides the estimate $||w||_{\mathcal{H}^2(\mathbb{R})} \lesssim ||P_{\pm k^*}[u_{\pm k^*}]||^2_{\mathcal{H}^2(\mathbb{R})}$. □

**Remark 12.** *We can simplify equation 5.1 a bit more if we take into account that $\mathcal{P}_{\pm k^*}$ and $\mathcal{L}_{\pm k^*}$ commute in the domain of $\mathcal{L}_{\pm k^*}$, by the definition of the operator $\mathcal{P}_{\pm k^*}$ we obtain*

$$\lambda_{\pm k^*}\mathcal{P}_{\pm k^*}[u_{\pm k^*}] = \mathcal{P}_{\pm k^*}\mathcal{N}_{\pm k^*}[\mathcal{P}_{k^*}u] \pm ik^*(\bar{d}-d)\mathcal{P}_{\pm k^*}u_{\pm k^*} \quad (5.4)$$

*which is finite dimensional, both in range and domain.*

5.1. **The finite dimensional system and its symmetries.** Given that the bifurcation equation in (5.4) is finite dimensional, we would like to use classical notion of continuity and differentiability to find a bifurcation. In order to do that we will need to know a few more properties of the eigenvalues $\lambda_{\pm k^*}$ and their associated eigenspaces:

**Proposition 10.** *The eigenvalues $\lambda_{\pm k^*}$ are differentiable, namely, $\epsilon \mapsto \lambda_{\pm k^*}(\epsilon) \in \mathscr{C}^1(\mathcal{I}; \mathbb{C})$. Furthermore, we can choose a basis $\{v_{\pm k^*}(\epsilon)\}$ associated to the eigenvalues $\lambda_{\pm k^*}(\epsilon)$ of the operators $\mathcal{L}_{\pm k^*}$ so that $\{v_{\pm k^*}(\epsilon)\} \in \mathscr{C}^1(\mathcal{I}; Range(\mathcal{P}_{\pm k^*}))$.*

*Proof.* This result follows from [CR73, section 1] in the complex case, in which it is proved that simple eigenvalues depend on the associated operator analiticaly (taking the operator norm). Since, as was proved in corollary 9, $(\epsilon, u) \mapsto \mathcal{L}_k^{\bar{d}}(\epsilon, u) \in \mathscr{C}^1(\mathcal{I} \times V; \mathcal{L}^2(\mathbb{R};\mathbb{R}^n))$, where $V \subset \mathcal{H}^2(\Omega;\mathbb{R}^n)$, the eigenvalues are $\mathscr{C}^1(\mathcal{I}; \mathbb{C})$ (shrink the interval $\mathcal{I}$ if necessary). The same article also proves that $v_{\pm k^*}(\epsilon) \in \mathscr{C}^1(\mathcal{I}; Range(\mathcal{P}_{\pm k^*}))$. □

Therefore, we can look for solutions $\mathcal{P}_{k^*}[u_{k^*}] = \zeta v_{k^*}(\epsilon)$; by the previous discussion, we can choose the eigenfunctions so that $v_{-k^*}(\epsilon) = \overline{v_{k^*}(\epsilon)}$. Applying this equation (5.3a), we obtain

$$\zeta\lambda_{k^*}(\epsilon)v_{k^*}(\epsilon) = f_{k^*}(\zeta v_{k^*}(\epsilon), \overline{\zeta}v_{-k^*}(\epsilon))v_{k^*}(\epsilon) + ik^*(\bar{d}-d)v_{k^*}(\epsilon) \quad (5.5a)$$

$$\overline{\zeta}\lambda_{-k^*}(\epsilon)v_{-k^*}(\epsilon) = f_{-k^*}(\zeta v_{k^*}(\epsilon), \overline{\zeta}v_{-k^*}(\epsilon))v_{-k^*}(\epsilon) - ik^*(\bar{d}-d)\overline{\zeta}v_{-k^*}(\epsilon) \quad (5.5b)$$

where $f_{\pm k^*}(\zeta v_{k^*}(\epsilon), \overline{\zeta}v_{-k^*}(\epsilon)) = \mathcal{P}_{\pm k^*}\mathcal{N}_{\pm k^*}[u]$. We focus on studying the action of the symmetries on $f_{\pm k^*}$. Initially, notice that the translational symmetry applies, so that we can state the next

**Proposition 11.** $(\zeta v_{k^*}(\epsilon), \overline{\zeta}v_{-k^*}(\epsilon))$ *solves* (5.5a) (5.5b) *if and only if* $(e^{i\theta}\zeta v_{k^*}(\epsilon), \overline{e^{i\theta}\zeta}v_{-k^*}(\epsilon))$ *solves it as well.*

It allow us to choose $\theta$ so that $e^{i\theta}\zeta$ is always real valued. In addition to that, $\mathcal{L}$ is a real operator and u is a real valued function taking real values,

$$\mathcal{L}_{-k}[u] = \frac{1}{2\pi}\int_0^{2\pi} e^{iky}\mathcal{L}[u(x,y)]dy = \overline{\frac{1}{2\pi}\int_0^{2\pi} e^{-iky}\mathcal{L}[u(x,y)]dy} = \overline{\mathcal{L}_k[u]}.$$

Applying the same reasoning to the operator $\mathscr{L}_{-k}[\cdot] - \lambda I$ one can see that $\overline{\sigma(\mathscr{L}_{-k})} = \sigma(\mathscr{L}_k)$. Thus, the eigenvalues for $\mathscr{L}_k$, $k \neq 0$, are forced to appear simultaneously. This symmetry in inherited by the spatial symmetry in the system. Further, from the symmetry $\overline{\mathscr{N}_{k^*}} = \mathscr{N}_{-k^*}$ we know that (5.5a) is the complex conjugate of (5.5b). Therefore, system of equations (1.11) is equivalent to

$$x\lambda_{k^*}(\epsilon)v_{k^*}(\epsilon) = f_{k^*}(xv_{k^*}(\epsilon), xv_{-k^*}(\epsilon))v_{k^*}(\epsilon)(\epsilon) + ik(\bar{d} - d)\zeta x v_{k^*}(\epsilon) \tag{5.6}$$

where $x \in \mathbb{R}$. Therefore, it suffices to solve

$$x\lambda_{k^*} = f_{k^*}(xv_{k^*}, xv_{-k^*}) + ik^*(\bar{d} - d)x \tag{5.7}$$

Depending on the underlying symmetry we take into account, this equation features a different behavior; however, we must begin by proving that $(\epsilon, x) \mapsto f_{k^*}(\epsilon, x)$ is differentiable in the classical sense:

**Proposition 12.** $(\epsilon, x) \mapsto f_{k^*}(\epsilon, x)$ is $\mathscr{C}^1(\mathcal{I} \times V; \mathbb{C})$, where $V \subset \mathbb{R}$ is a neighborhood of the origin.

*Proof.* By te Hanh- Banach theorem, we can choose a dual element $h$ so that $h(v_{k^*}(0)) = 1$. As $h$ is a bounded linear functional, $h \in \mathscr{C}^\infty(\mathscr{H}^2(\mathbb{R}; \mathbb{R}^n); \mathbb{C})$. Applying h on both sides of the system (5.6) and shrinking the interval $\mathcal{I}$ if necessary so that $h(v_{k^*}(\epsilon) \neq 0$, $\forall \epsilon \in \mathcal{I}$ we arrive at (5.7). □

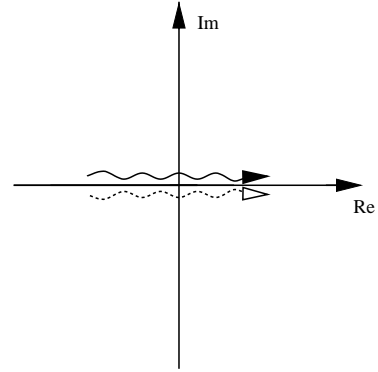

FIGURE 2. Eigenvalues collapsing at the real line: O(2) symmetry implies double eigenvalues

The next point we want to stress is a spectral property in the O(2)-symmetric setting:

**Proposition 13.** *The bifurcating eigenvalues are real in the O(2) symmetric case*

In the O(2)-symmetric setting this implies:

**Proposition 14.** *In the O(2) symmetric setting the mapping $f_{k^*}$ as given by equation (5.7) is real valued.*

In order to prove these propositions we must understand a bit more the relation between the symmetries and the projections $\mathcal{P}_k$ we have been using.

**Claim 4.** *The following properties hold: $\mathcal{P}_k \Gamma(u)] = R\mathcal{P}_{-k}[u]$, $R\mathscr{L}_{-k} = \mathscr{L}_k R$ and $R\mathcal{P}_{-k}[\cdot] = \mathcal{P}_k R[\cdot]$.*

*Proof.* The first property follows by standard calculations. Since $\Gamma\mathscr{L} = \mathscr{L}\Gamma$, we have that

$$\mathcal{P}_k[\Gamma\mathscr{L}] = (\mathcal{P}_k \Gamma)\mathscr{L} = R\mathcal{P}_{-k}\mathscr{L} = R(\mathscr{L}_{-k}\mathcal{P}_{-k})$$

On the other hand, $\mathcal{P}_k[\mathscr{L}\Gamma] = \mathscr{L}_k(\mathcal{P}_k\Gamma) = \mathscr{L}_k R\mathcal{P}_{-k}$. We arrive at the second property, since $\mathcal{P}_k$ is surjective. We obtain the last property by

$$R\mathcal{P}_k[f] = R\frac{1}{2\pi i}\int_\gamma (\lambda I - \mathscr{L}_k)^{-1} f d\lambda = \frac{1}{2\pi i}\int_\gamma (\lambda I - \mathscr{L}_{-k})^{-1} Rf d\lambda = P_{-k}[Rf].$$

□

**Proposition 15.** *The finite dimensional reduction inherits the translation symmetry in both cases; in particular, it inherits the reflection symmetry in the O(2) case.*

*Proof.* Analogue to the proofs on **section 4**, since we have uniqueness of the contraction mapping and the neighborhoods are still the same upon symmetry action. □

*Proof of proposition 13.* Let $(\lambda_k, u_k)$ be a pair eigenvalue/eigenvector of $\mathscr{L}_k$. Then, $u_k(x)e^{iky}$ is an eigenvalue of $\mathscr{L}$. By the reflection symmetry, we know that $\Gamma[u_k e^{iky}]$ is also an eigenvalue for $\mathscr{L}$ associated to $\lambda_k$. By property 4, we obtain

$$\Gamma(\lambda_k[u_k(x)e^{iky}]) = \lambda_k e^{-iky} R[u_k(x)] = \Gamma\mathscr{L}[u_k(x)e^{iky}] = \mathscr{L}[\Gamma[u_k(x)e^{iky}]] =$$
$$= \mathscr{L}[e^{-iky} R u_k(x)] = e^{-iky} \mathscr{L}_{-k}[R u_k(x)]$$

As $\lambda_{-k}$ is the only bifurcating eigenvalue of $\mathscr{L}_{-k}$, we get that $\lambda_k = \lambda_{-k}$. By the precedent discussion, we also know that $\bar{\lambda}_k = \lambda_{-k}$. Therefore $\bar{\lambda}_k = \lambda_k \Rightarrow \lambda_k \in \mathbb{R}$. $\square$

*Proof of proposition 14.* Equation 5.7 in this case reduces to $x\lambda_{k^*} = f_{k^*}(xv_{k^*}, xv_{-k^*})$. By hypotheses (H1), we can use the splitting of the spectrum to define projectors. Therefore, the equality

$$(\lambda_{k^*} I - \mathscr{L}_{k^*})^{-1} \mathscr{N}_{k^*}[\zeta v_{k^*}(\epsilon), \overline{\zeta} v_{-k^*}(\epsilon)] = \frac{1}{2\pi i} \int_\gamma \mathcal{P}_{k^*}[(\lambda I - \mathscr{L})^{-1} \mathscr{N}[\zeta v_{k^*}(\epsilon), \overline{\zeta} v_{-k^*}(\epsilon)]]$$

holds, by the definition of $\mathcal{P}_k[\cdot]$, since the right hand side is defined. As discussed in the proof of **proposition 13**, we know that $e^{\pm ik^* y} v_{k^*}(x)$ span a basis to the range of $\mathcal{P}_{k^*}$. Since we want real valued solutions, we define

$$V_1 := \frac{e^{iky} P_k(u_k) + e^{-iky} P_{-k}(u_{-k})}{2} \quad, \quad V_2 := \frac{e^{iky} P_k(u_k) - e^{-iky} P_{-k}(u_{-k})}{2i}.$$

as a basis to the eigenspace of the operator $\mathscr{L}$ associated to the eigenvalue $\lambda_{k^*} (= \lambda_{-k^*})$. Thus,

$$\frac{1}{2\pi i} \int_\gamma (\mu I - \mathscr{L})^{-1} \mathscr{N}(u) d\mu = f(u) V_1 + g(u) V_2,$$

where f and g are real (since only real valued functions are to be taken as solutions of 1.11). Now, applying $\mathcal{P}_k$ on both sides and using the lemma and using the definition of $f_{k^*}$ we get

$$P_k \mathscr{N}_k(u) d\mu = f_{k^*}(\zeta v_{k^*}(\epsilon), \overline{\zeta} v_{-k^*}(\epsilon)) = \frac{f(\zeta v_{k^*}(\epsilon), \overline{\zeta} v_{-k^*}(\epsilon))}{2} P_k + \frac{g(\zeta v_{k^*}(\epsilon), \overline{\zeta} v_{-k^*}(\epsilon))}{2i} P_k.$$

Similarly, $P_{-k} \mathscr{N}_{-k}(u) d\mu = \frac{f(u)}{2} P_{-k} - \frac{g(u)}{2i} P_{-k}$. Conjugating and using the symmetries we have discussed just below **proposition 11**, we end up with this system:

$$\lambda P_k[u_k] = \frac{f(\zeta v_{k^*}(\epsilon), \overline{\zeta} v_{-k^*}(\epsilon)) - ig(\zeta v_{k^*}(\epsilon), \overline{\zeta} v_{-k^*}(\epsilon))}{2} P_k[u_k]$$
$$\lambda P_{-k}[u_{-k}] = \frac{f(\zeta v_{k^*}(\epsilon), \overline{\zeta} v_{-k^*}(\epsilon)) + ig(\zeta v_{k^*}(\epsilon), \overline{\zeta} v_{-k^*}(\epsilon))}{2} P_{-k}[u_{-k}].$$

Applying R (which is invertible) to the second equation gives

$$\lambda R[P_{-k}(u_{-k})] = \frac{f(\zeta v_{k^*}(\epsilon), \overline{\zeta} v_{-k^*}(\epsilon)) + ig(\zeta v_{k^*}(\epsilon), \overline{\zeta} v_{-k^*}(\epsilon))}{2} R[P_{-k}(u_{-k})] \Rightarrow$$
$$\Rightarrow \lambda[P_k(u_k)] = \frac{f(\zeta v_{k^*}(\epsilon), \overline{\zeta} v_{-k^*}(\epsilon)) + ig(\zeta v_{k^*}(\epsilon), \overline{\zeta} v_{-k^*}(\epsilon))}{2} P_k(u_k).$$

Therefore,

$$\lambda P_k[u_k] = \frac{f((\zeta v_{k^*}(\epsilon), \overline{\zeta}v_{-k^*}(\epsilon)) - ig(\zeta v_{k^*}(\epsilon), \overline{\zeta}v_{-k^*}(\epsilon))}{2} P_k[u_k]$$
$$\lambda P_k[u_k] = \frac{f(\zeta v_{k^*}(\epsilon), \overline{\zeta}v_{-k^*}(\epsilon)) + ig(\zeta v_{k^*}(\epsilon), \overline{\zeta}v_{-k^*}(\epsilon))}{2} P_k[u_k] \quad \Rightarrow \quad g(\zeta v_{k^*}(\epsilon), \overline{\zeta}v_{-k^*}(\epsilon)) = 0.$$

proving that $f_{k^*}(\zeta v_{k^*}(\epsilon), \overline{\zeta}v_{-k^*}(\epsilon))$ is real valued in the O(2) setting, as we sought. $\square$

## 6. The bifurcation analysis

Equation (5.7) has a different nature depending on the symmetry we take into account. We treat the two cases separatedly then.

### 6.1. The $\mathcal{O}(2)$ case: proof of theorem 1.13.
Consider the bifurcation equations (5.7). Let $f(x) = f_{k^*}(xv_{k^*}, xv_{-k^*})$. Define the equation $F(\epsilon, x) = \lambda(\epsilon)x - f(x)$. Clearly, $F(\epsilon, 0) = 0$, for $\epsilon \in \mathcal{I}$. We want to apply Crandal-Rabinowitz theorem on bifurcation on a simple eigenvalue ([CR71]). Notice that $D_x F\big|_{(\epsilon,x)=(0,0)} = \lambda(0) - D_x f(0,0) = \lambda(0) = 0$. Thus, $\hat{1} \in ker D_x F$, $Ran[D_x f] = 0$. On the other hand,

$$D_\epsilon D_x F\big|_{(\epsilon,x)=(0,0)} = \lambda'(\epsilon)\big|_{(\epsilon,x)=(0,0)} + D_\epsilon D_x f\big|_{(\epsilon,x)=(0,0)} = \lambda'(\epsilon) \neq 0.$$

As $D_\epsilon D_x F\big|_{(\epsilon,x)=(0,0)}[\hat{1}] \notin Ran[D_x F|_{(\epsilon,x)=(0,0)}]$, we can apply the bifurcation on simple eigenvalue theorem to get a nontrivial solution parametrized on a small open interval $\mathcal{I}$ around 0 such that

$$I \ni s \mapsto F[\epsilon(s), x(s)] = 0.$$

### 6.2. The $\mathcal{SO}(2)$ case: proof of theorem 1.12.
In this case there is no reflection symmetry, so we don't have $\lambda$ real in general. Further, the last term in the right hand side of (5.7) is not necessarily zero. The problem can be written as

$$\lambda x = f(x, \bar{x}) + i(\bar{d} - d)k^* x \quad \text{or} \quad \begin{cases} Re[\lambda]x = Re[f(x)] \\ Im[\lambda]x = Im[f(x)] + (\bar{d} - d)k^* x \end{cases}.$$

Define then

$$G[\epsilon, d, x] = \begin{bmatrix} Re[\lambda(\epsilon)]x - Re[f(x)] \\ Im[\lambda(\epsilon)]x - Im[f(x)] + (\bar{d} - d)k^* x \end{bmatrix}$$

We know that $G[\epsilon, d, 0] = 0$. But now we have two parameters, so we cannot apply the theorem for bifurcation on a simple eigenvalue as stated in [CR71]. Notice, however, that we can apply it on the first equation which, once solved, can be plugged to the second equation, reducing the problem to the 1-D case. In detail: as in the previous case, there exists an open interval $\mathcal{I}$ around 0 such that $\mathcal{I} \ni s \mapsto G_1[\epsilon(s), d, x(s)] = G_1[\epsilon(s), x(s)] = 0$. Plugging in the second equation, we obtain $I \ni s \mapsto G_2[\epsilon(s), d, x(s)]$, which we know to satisfy $G_2[\epsilon(0), d, x(0)] = 0$. So, define $f(s, d) = G_2[\epsilon(s), d, x(s)]$. We know that $f(0, d) = 0$. Also,

$$\frac{\partial f}{\partial s}\bigg|_{(s,d)=(0,0)} = Im[\lambda'(0)]x(0) - Im[\lambda(0)]x'(0) - D_\epsilon Im f(0,0)\epsilon'(0) - D_x Im f(0,0)x'(0) + (\bar{d} - d)k^* x'(0).$$

Where, $x'(0)$ is tangent to the kernel (i.e., $x'(0) \neq 0$). As $x(0) = 0, \lambda(0) = 0$, we obtain

$$\frac{\partial f}{\partial s}\bigg|_{(s,d)=(0,0)} = +(\bar{d} - d)k^* x'(0) \neq 0.$$

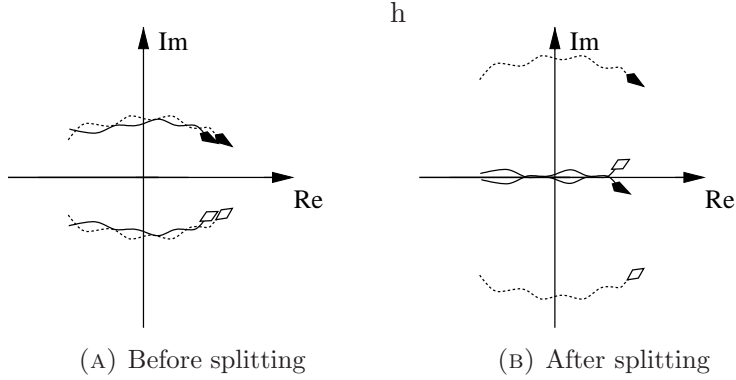

(A) Before splitting  (B) After splitting

FIGURE 3. O(2)-Hopf bifurcation: double eigenvalues are separated by the Galilean change of frame in y-direction.

Thus there exists a $\mathscr{C}$ mapping $s \mapsto d(s)$ s.t. $f(s, d(s)) = G_2[\epsilon(s), d(s), x(s)] = 0$ in a neighborhood of 0, and we are done.

### 6.3. Proof of the corollary 1: the O(2) Hopf-bifurcation.
Initially, notice that the O(2) symmetry creates double eigenvalues: if $v(\epsilon)$ is an eigenvectur associated with $\lambda(\epsilon)$, $\Gamma[v(\epsilon)]$ is also an eigenvalue. By assumption, the crossing eigenvalues are of the transversal type, i.e., associated with nonzero modes ($k^* \neq 0$). Applying the projector $\Pi_{k^*}$ on $\mathscr{L}^{\bar{d}}[v(\epsilon)] = \lambda(\epsilon)v(\epsilon)$ and $\mathscr{L}^{\bar{d}}[\Gamma v(\epsilon)] = \lambda(\epsilon)\Gamma v(\epsilon)$ and using some of properties given by claim 4, we obtain

$$\mathscr{L}^{\bar{d}}_{k^*}[v(\epsilon)] = \lambda(\epsilon)v(\epsilon) \quad , \quad \mathscr{L}^{\bar{d}}_{-k^*}[\Gamma[v(\epsilon)]] = \lambda(\epsilon)\Gamma[v(\epsilon)]$$

Hence these eigenvalues are associated with different modes ($\pm k^*$). Upon a Galilean change of frame in y-direction (see section 1.4), one of the later eigenvalues is mapped to the origin, while other is mapped to $2\lambda(0)$ (figures 3a and 3b). Notice that, in a neighborhood of the origin we have the same picture as before: a pair of eigenvalues crossing the origin. As the result on appendix D just uses projections as integrals of contours of the origin, we can ignore the rest of the spectrum using the same technique we applied in section 5, ending up in the same context as that treated by theorem 1.13, since the O(2)-symmetry has been broken (see remark 5).

Concerning the solution v(s), as the later is steady in the moving frame $(x - ct, y - d(s)t)$ and the domain $\Omega$ is $2\pi$ periodic in y-direction, one can see that the solution in the frame $(x - ct, y)$ is $\dfrac{2\pi}{d(s)}$ periodic for s sufficiently small[10], asserting the existence of nontrivial periodic orbits.

## APPENDIX A. THE PARAMETRIC CONTRACTION MAPPING THEOREM

**Theorem A.1.** *Let $\mathscr{T} : V \times B \mapsto B$, where B is a Banach space, V is an open set in a topological space such that $||\mathscr{T}(\lambda, x) - \mathscr{T}(\lambda, y)|| \leq \theta ||x - y||, \theta < 1$, for every $\lambda \in V$, $x, y \in B$. Then, we know that for every $\lambda \in V$ and $z \in B$, there exists a solution $x = x(\lambda, z)$ to the problem*

$$x = z + T(\lambda, x)$$

*Furthermore, if $T(\cdot, y)$ is a $\mathscr{C}^k$ mapping on both entries then $V \ni \lambda \mapsto x(\lambda)$ is also $\mathscr{C}^k$.*

*Proof.* see [CH82] and their proof of the implicit function theorem, or [Ham82, section I]. □

---
[10]Notice that, as $\bar{d} = d(0) \neq 0$, $d(s) \neq 0$ in in a neighborhood of 0 in $\mathcal{J}$.

**Corollary 6** (Bounds on implicit functions). *Let $T : X \times Y \to X$, continuous on both entries. Assume the $T(x,y)$ has a unique fixed point $x$ for every $y$ ( i.e., $y \mapsto x(y)$ is well defined) and that $|T(x,y)|_X \leq C(|x|_X^k + |y|_Y^k)$. Then, $|x(y)| \lesssim |y|^k$*

*Proof.* As T is continuous, $y \mapsto x(y)$ is also continuous. Thus, for $\delta$ sufficiently small, we can take $C||x(y)||_X^{k-1} < \dfrac{1}{2}$. Then,

$$x(y) = T[x(y), y] \Rightarrow ||x(y)||_X \leq C(||x(y)||_X^k + ||y||_Y^k) \leq \frac{1}{2}||x(y)||_X + C||y||_Y^k \Rightarrow ||x(y)||_X \leq 2C||y||_Y^k$$

$\square$

**Corollary 7** (Perturbation of the contraction by small bounded operators). *Let $T : X \to X$ be a contraction with parameter $\theta < 1$. Then, for all sufficiently small bounded operator $P$, $T + P$ is also a contraction.*

*Proof.* Being a contraction, T has the following properties:there exists a ball B such that $T(B) \subset B$; T has Lipschitz constant $\theta < 1$. It is not hard to prove that the ball B can be taken with radius r so that $T(B[r]) \subset B[r - \epsilon]$, for an small $0 < \epsilon < r$. We prove that $(T + P)(B[r]) \subset B[r]$; let $x \in B[r]$. Let $||P|| \leq \dfrac{\epsilon}{r} \wedge \dfrac{1 - \theta}{2}$. Then,

$$||(T + P)(x)|| = ||T(x) + P(x)|| \leq (||T(x)|| + ||P||||x||) \leq (r - \epsilon) + ||P||r \leq \leq (r - \epsilon) + \delta_1 r = r$$

It remains to be proved that T+ P is a contraction. But it is a consequence of P's linearity, as one can see from

$$||(T + P)[u - v]|| \leq ||T[u - v]|| + ||P[u - v]|| \leq (\theta + ||P||)||u - v|| \leq \frac{\theta + 1}{2}||(T + P)[u - v]||.$$

$\square$

APPENDIX B. USEFUL FACTS USED IN THE PROOF

B.1. **A few properties of the mapping $\Pi_k[\cdot]$.**

**Proposition 16.** *Let $u \in \mathscr{H}^s(\Omega; \mathbb{R}^n)$ and $\alpha$ be a multi-index such that $|\alpha| \leq s$. Then $\partial^\alpha(\tilde{\Pi}_k[u]) = (ik)^{\alpha_y}\tilde{\Pi}_k[\partial^\alpha u]$ in the weak sense.*

*Proof.* Given $\phi \in \mathscr{C}^\infty(\mathbb{R}; \mathbb{R}^n)$,

$$\begin{aligned}(-1)^\alpha \int_\Omega \tilde{\Pi}_k[u](x,y) \cdot \partial^\alpha \phi(x,y) dx &= (-1)^\alpha \int_\Omega \left\{\frac{1}{2\pi} \int_\mathbb{T} u(x,\xi)e^{-ik\xi}d\xi\right\} \cdot e^{iky}\partial^\alpha \phi(x,y) dxdy = \\ &= (-1)^\alpha \int_\mathbb{T} \left\{\frac{1}{2\pi} \int_\Omega u(x,\xi) \cdot e^{-ik\xi}\partial^\alpha \phi(x,y) dxd\xi\right\} e^{iky}dy\end{aligned}$$

As $e^{-ik\xi}\partial^\alpha \phi(x,y) = \partial_x^{\alpha_x}[e^{-ik\xi}\partial_y^{\alpha_y}\phi(x,y)]$, the last integral is equal to

$$(-1)^\alpha \int_{\mathbb{T}} \left\{ \frac{1}{2\pi} \int_\Omega \partial_x^{\alpha_x} u(x,\xi) \cdot e^{-ik\xi} \partial_y^{\alpha_y} \phi(x,y) dx d\xi \right\} e^{iky} dy =$$

$$= \frac{(-1)^{\alpha_y}}{2\pi} \int_{\mathbb{T}} \partial_y^{\alpha_y} \left\{ \int_\Omega \partial_x^{\alpha_x} u(x,\xi) \cdot e^{-ik\xi} \phi(x,y) dx d\xi \right\} e^{iky} dy =$$

$$\stackrel{Int.\,by.\,parts}{=} \frac{(ik)^{\alpha_y}}{2\pi} \int_{\mathbb{T}} \left\{ \int_\Omega \partial_x^{\alpha_x} u(x,\xi) \cdot e^{-ik\xi} \phi(x,y) dx d\xi \right\} e^{iky} dy =$$

$$= (ik)^{\alpha_y} \int_\Omega \tilde{\Pi}_k[\partial_x^{\alpha_x} u](x,y) \cdot \phi(x,y) dx dy$$

□

**Theorem B.1.** $\tilde{\Pi}_k \circ \mathscr{L}[u](x,y) = \mathscr{L}_k \circ \Pi_k[u](x)e^{ik \cdot y}$ a.e. since, for any given $\phi \in \mathscr{C}_c^\infty(\mathbb{R} \times \mathbb{T}; \mathbb{R}^n)$,

*Proof.* Let $\alpha = (\alpha_x, \alpha_\xi)$ be a multi-index.

$$\int_\Omega \{\tilde{\Pi}_k \circ \mathscr{L}[u](x,y)\} \cdot \phi(x,y) dx dy = \int_\Omega \left\{ \int_{\mathbb{T}} \mathscr{L}[u](x,\xi) e^{-ik\xi} d\xi \right\} \cdot e^{iky} \phi(x,y) dx dy =$$

$$= \int_\Omega \left\{ \int_{\mathbb{T}} (\sum_\alpha \partial_x^{\alpha_x} \partial_\xi^{\alpha_\xi} [A_\alpha(x) u(x,\xi)] \cdot (e^{-ik\xi} \phi(x,y)) d\xi \right\} e^{iky} dx dy =$$

$$\stackrel{Fubini}{=} \int_{\mathbb{T}} \left\{ \frac{1}{2\pi} \int_\Omega \sum_\alpha \partial_x^{\alpha_x} \partial_\xi^{\alpha_\xi} [A_\alpha(x) u(x,\xi)] \cdot (e^{-ik\xi} \phi(x,y)) dx d\xi \right\} e^{iky} dy =$$

$$\int_{\mathbb{T}} \left\{ \frac{(-1)^{\alpha_\xi}}{2\pi} \int_\Omega \sum_\alpha \partial_x^{\alpha_x} [A_\alpha(x) u(x,\xi)] \cdot \partial_\xi^{\alpha_\xi}(e^{-ik\xi} \phi(x,y)) dx d\xi \right\} e^{iky} dy =$$

$$\int_{\mathbb{T}} \left\{ \frac{(ik)^{\alpha_\xi}}{2\pi} \int_\Omega \sum_\alpha \partial_x^{\alpha_x} [A_\alpha(x) u(x,\xi)] \cdot (e^{-ik\xi} \phi(x,y)) dx d\xi \right\} e^{iky} dy =$$

Above, we used that $(x,\xi) \mapsto A_\alpha(x) e^{-ik\xi} \phi(x,y) \in \mathscr{C}_c^\infty(\mathbb{R} \times \mathbb{T}; \mathbb{R}^n)$ for every $y \in \mathbb{T}$, and that $A(x,\xi) = A(x)$. Using the result of lemma 1 - which says that $\partial_x \alpha_x \Pi_k[u] = \Pi_k[\partial_x \alpha_x u]$ and is is independent of the present proof - we obtain

$$\int_{\mathbb{T}} \left\{ \frac{(ik)^{\alpha_\xi}}{2\pi} \int_\Omega \sum_\alpha \partial_x^{\alpha_x} [A_\alpha(x) u(x,\xi)] \cdot (e^{-ik\xi} \phi(x,y)) dx d\xi \right\} e^{iky} dy =$$

$$\int_{\mathbb{T}} \left\{ \frac{(ik)^{\alpha_\xi}}{2\pi} \sum_\alpha \partial_x^{\alpha_x} \int_\Omega [A_\alpha(x) u(x,\xi)] \cdot (e^{-ik\xi} \phi(x,y)) dx d\xi \right\} e^{iky} dy$$

Another application of Fubini gives

$$(ik)^{\alpha_\xi} \int_\Omega \sum_\alpha \partial_x^{\alpha_x} A_\alpha(x) \underbrace{\left[\frac{1}{2\pi} \int_{\mathbb{T}} u(x,\xi) e^{-ik\xi} d\xi\right]}_{u_k(x)} \cdot \phi(x,y) e^{iky} dx dy = \int_\Omega \mathscr{L}_k \left\{ \Pi_k[u](x) e^{ik \cdot y} \right\} \cdot \phi(x,y) dx dy$$

Summarizing,

$$\int_\Omega \{\tilde{\Pi}_k \circ \mathscr{L}[u](x,y)\} \cdot \phi(x,y) dx dy = \int_\Omega \mathscr{L}_k \circ \Pi_k[u](x) e^{ik \cdot y} \cdot \phi(x,y) dx dy$$

$\square$

APPENDIX C. DIFFERENTIABILITY WITH RESPECT TO PARAMETERS

We summarize a few of the Frechet-differentiability results we have used from the perspective of substitution operators. This appendix, besides making the article self contained, intends to make the use of these results more clear (see [Cra77] and references therein; we also recommend the excellent article [Ham82]).

C.1. **Substitution operators.** It is well know that even $\mathscr{C}^\infty(\mathbb{R};\mathbb{R})$ mappings may not be continuous when composed with Sobolev functions (see [Cra77]). We devote this section to showing that the mappings we have been using (as $\mathscr{L}$, $\mathscr{N}$) and their variants have differentiable properties in the appropriate spaces. Essentially, it concerns to understanding how the space $\mathscr{C}^k(\mathbb{R}^n,\mathbb{R}^n)$ acts on the space $\mathscr{H}^2(\Omega;\mathbb{R}^n)$; as we will see soon, the embedding theorem of **section 2.3** is one of the main ingredients.

**Lemma 4.** *Let $x \mapsto f(x)$ be so that $f \in \mathscr{C}^k(\mathbb{R}^n;\mathbb{R}^m)$ and $K \subset \mathbb{R}^n$ be a compact set. Then, for every fixed $\eta > 0$ there exists a $\delta > 0$ such that*

$$|f(x+h) - \sum_{\alpha \leq k} \frac{1}{\alpha!} \partial_x^\alpha f(x,) h^\alpha| \leq \eta |h|^k \quad , \quad |h| \leq \delta, \forall x \in K$$

*Proof.* Define the mapping $g(x,h) := \frac{1}{|h|^k}(f(x+h) - \sum_{\alpha \leq k} \frac{1}{\alpha!} \partial_x^\alpha f(x,) h^\alpha)$ if $|h| \neq 0$, $g(x,0) = 0$ otherwise. As g is continuous, then g is uniformly continuous over compact sets. Therefore, given an $\epsilon > 0$ and a compact set $K \times [0,1]$, there exists a $c > 0$ such that $|g(x,h) - g(\tilde{x},\tilde{h})| \leq \eta$ whenever $|x - \tilde{x}| - |h - \tilde{h}| \leq \delta$. Recalling that $g(x,0,y) = 0$ we get that $|g(x,h,y)| \leq \epsilon$ for all $|h| \leq \delta$. $\square$

We begin with the following result:

**Lemma 5.** *Let $\mathcal{R}^{(i)}(\cdot,\cdot)$ as defined by equation (1.3b); then, $\forall i \in \{1,2\}$,*

$$\mathcal{R}^{(i)} \in \mathscr{C}^1(\mathcal{I} \times V; \mathscr{H}^1(\Omega;\mathbb{R}^n)),$$

*where V is a neighborhood of zero in $\mathscr{H}^2(\Omega;\mathbb{R}^n)$*

*Proof.* Fix $(\tilde{\epsilon}, u)$ in a neighborhood of $(0,0) \in \mathcal{I} \times \mathscr{H}^2(\Omega;\mathbb{R}^n)$. Without loss of generality, we denote this neighborhood by $\mathcal{I} \times V$ (shrink the interval $\mathcal{I}$ if necessary), $0 \in V \subset \mathscr{H}^2(\Omega;\mathbb{R}^n)$. Given $(\epsilon, v) \in \mathbb{R} \times \mathscr{H}^2(\Omega;\mathbb{R}^n)$, we need estimate

$$\mathcal{R}^{(i)}(\epsilon + \tilde{\epsilon}, u + v) - \mathcal{R}^{(i)}(\tilde{\epsilon}, u) - \epsilon \partial_\epsilon \mathcal{R}^{(i)}(\tilde{\epsilon}, u) - v \partial_u \mathcal{R}^{(i)}(\tilde{\epsilon}, u) =: \Theta(\epsilon, v)$$

which, in this case, can be written as

$$\Theta(\epsilon, v) := \int_0^1 \frac{d^2}{ds^2} \mathcal{R}^{(i)}(\epsilon + s\tilde{\epsilon}, u + sv)(1-s) ds.$$

We must show that $\Theta(\epsilon, v) = o(|\epsilon| + ||v||_{\mathscr{H}^2(\Omega)})$. In fact, it can be rewritten as

$$\Theta(\epsilon, v) = I\epsilon^2 + 2II\epsilon(v) + III(v,v),$$

where $I = \int_0^1 \partial_\epsilon^2 \mathcal{R}^{(i)}(\epsilon + s\tilde{\epsilon}, u + sv)(1-s)ds$, $v \mapsto II(v) = \int_0^1 \partial_\epsilon \partial_u \mathcal{R}^{(i)}(\epsilon + s\tilde{\epsilon}, u + sv)(1-s)ds(v)$ is a linear mapping and $(u, w) \mapsto III(u, w) = \int_0^1 \partial_u^2 \mathcal{R}^{(i)}(\epsilon + s\tilde{\epsilon}, u + sv)(1-s)ds(u, w)$ is a bilinear mapping. We will prove the result by proving that $\Theta(\epsilon, v) = O(|\epsilon|^2 + ||v||^2_{\mathscr{H}^2(\Omega)})$. As $f^{(1)}$, $f^{(2)}$ are $\mathscr{C}^2(\mathcal{I} \times \mathbb{R}^n; \mathbb{R}^n)$, an application on lemma 4, jointly with the embedding theorem of section 2.3, gives that I and II are in $\mathscr{L}^\infty(\Omega; \mathbb{R}^n)$, with bound uniform in $\mathcal{I} \times V$; therefore, $II(v)$ and $III(v,v)$ are $\mathscr{L}^2(\Omega; \mathbb{R}^n)$ integrable, i.e., there exists a constant C and a $\delta > 0$ such that $II \leq C\epsilon ||v||_{\mathscr{H}^2(\Omega)}$ and $III \leq C||v||^2_{\mathscr{H}^2(\Omega)}$. In order to estimate $I$, notice that its integrand is bounded by $C(|v + u|^2)$. Therefore, we can compute the $\mathscr{L}^2$ integral, and the result is achieved. The estimate for the weak derivative (which we know to exist, by claim 2) is similarly obtained. $\square$

**Corollary 8.** *Let $\mathcal{N}$ be defined by (1.3a); then*
$$\mathcal{N} \in \mathscr{C}^1(\mathcal{I} \times V; \mathscr{L}^2(\Omega; \mathbb{R}^n)),$$
*where $V \subset \mathscr{H}^2(\Omega; \mathbb{R}^n)$ is a neighborhood of the origin. Furthermore, the operator $\tilde{\mathcal{N}}_0$ as defined in $\tilde{\mathcal{N}}_0 \in \mathscr{C}^1(\mathcal{I} \times \mathscr{H}^2(\Omega; \mathbb{R}^n); \mathscr{H}^1(\Omega; \mathbb{R}^n))$.*

*Proof.* Use the definition of $\mathcal{N}$ and $\mathcal{R}^{(i)}$ for $i \in \{1, 2\}$ from (1.3a) and (1.3b) and the fact that the operator $\partial_x$ is $\mathscr{C}^\infty$ from $\mathscr{H}^1(\Omega; \mathbb{R}^n)$ to $\mathscr{L}^2(\Omega; \mathbb{R}^n)$. $\square$

We prove now results concerning to the operator $\mathscr{L}^{\bar{d}}$:

**Lemma 6.** *Let $\mathscr{L}^{\bar{d}}$ be defined as in (1.7). Then, there exists a neighborhood $\mathcal{I} \times V \ni (0, 0)$ in $\mathbb{R} \times \mathscr{H}^2(\Omega; \mathbb{R}^n)$ such that*
$$\mathscr{L}^{\bar{d}} \in \mathscr{C}^1(\mathcal{I} \times V; \mathscr{L}^2(\Omega; \mathbb{R}^n)).$$

*Proof.* Just notice that the matrices $A(,)$ and $\partial_\epsilon A(\cdot, \cdot)$ are locally bounded in the $\mathscr{L}^\infty$ norm. Therefore, for every $\epsilon \in \mathcal{I}$ and $v, w \in V \subset V$ sufficiently small, $A(\epsilon, w)v \in \mathscr{L}^2(\Omega; \mathbb{R}^n)$. Now, given $(\tilde{ep}, u) \in \mathcal{I} \times V$ and $(\epsilon, u) \in \mathbb{R} \times \mathscr{H}^2(\Omega; \mathbb{R}^n)$ we estimate
$$\mathscr{L}^{\bar{d}}(\epsilon + \tilde{\epsilon}, u + v) - \mathscr{L}^{\bar{d}}(\epsilon + \tilde{\epsilon}, u) - \epsilon \partial_\epsilon \mathscr{L}^{\bar{d}}(\tilde{\epsilon}, u) - \mathscr{L}^{\bar{d}}(\tilde{\epsilon}, v).$$
which we rewrite as $\{\mathscr{L}^{\bar{d}}(\epsilon + \tilde{\epsilon}, u) - \mathscr{L}^{\bar{d}}(\tilde{\epsilon}, u) - \epsilon \partial_\epsilon \mathscr{L}^{\bar{d}}(\tilde{\epsilon}, u)\} + \{\mathscr{L}^{\bar{d}}(\epsilon + \tilde{\epsilon}, v) - \mathscr{L}^{\bar{d}}(\tilde{\epsilon}, v)\} = I + II$. As $|I| = |A(\epsilon + \tilde{\epsilon}) - A(\tilde{\epsilon}) - \epsilon \partial_\epsilon A(\tilde{\epsilon})| \cdot |u|$ for an appropriate A, we have $\mathscr{L}^2(\Omega; \mathbb{R}^n)$ integrability and, overall, $||I|| = o(|\epsilon| + ||v||_{\mathscr{H}^2(\Omega)})$. On the other hand, $|II| \leq |A(\epsilon + \tilde{\epsilon}) - A(\epsilon)| \cdot |v|$; hence $||II||_{\mathscr{L}^2(\Omega; \mathbb{R}^n)} \leq ||A(\epsilon + \tilde{\epsilon}) - A(\epsilon)||_{\mathscr{L}^\infty(\Omega)} \cdot ||v||_{\mathscr{L}^2(\Omega; \mathbb{R}^n)}$. Uniform continuity of A over compact sets implies the result. $\square$

The next part of this appendix is concerned with the implications of differentiability as obtained upon composition with the projections $\Pi_j$. We start with

**Proposition 17.** $\Pi_j[\cdot] \in \mathscr{C}^\infty(\mathscr{H}^2(\Omega; \mathbb{R}^n); \mathscr{H}^2(\mathbb{R}; \mathbb{C}^n))$, *and the bound C on $||\Pi_j[u]||_{\mathscr{H}^2(\mathbb{R})} \leq C||u||_{\mathscr{H}^2(\Omega)}$ is uniform in $j \in \mathbb{Z}$. Moreover, $S[\cdot] \in \mathscr{C}^\infty(\mathscr{H}^2(\Omega; \mathbb{R}^n); M_{(0, \pm k^*)})$.*

*Proof.* Use linearity and lemma 1, namely, $\partial_x^\alpha \Pi_j[u] = \Pi_j[\partial_x^\alpha u]$ a.e.. $\square$

As a consequence of the composition of the operators, we obtain a

**Corollary 9.** *Let $\mathscr{L}_j^{\bar{d}} = \Pi_j \circ \mathscr{L} : \mathcal{I} \times \mathscr{H}^2(\Omega; \mathbb{R}^n) \to \mathscr{L}^2(\Omega; \mathcal{C}^n)$. Then,*
$$\mathscr{L}_j^{\bar{d}} \in \mathscr{C}^1(\mathcal{I} \times V; \mathscr{L}^2(\Omega; \mathcal{C}^n)),$$
*where $V \subset \mathscr{H}^2(\Omega; \mathbb{R}^n)$ can be chosen uniformly on $j \in \mathbb{Z}$.*

*Proof.* $\Pi_j[\cdot]$ is independent of $\epsilon$ and $j \in \mathbb{Z}$, thus we have uniform bounds. The result follows by differentiability of each term of the composed function $\mathscr{L}_j^{\bar{d}}$, $\mathscr{L}^{\bar{d}}$ and $\Pi_j$. $\square$

The rest of this appendix is dedicated to the following proof:

**Theorem C.1.** *Let the mapping $T : \mathcal{I} \times M_{(0,\pm k^*)} \to M_{(0,\pm k^*)}$ be defined as in theorem 3.4. Then, there exists a neighborhood $\mathcal{I} \times W \ni (0,0)$ such that*

$$T \in \mathscr{C}^1(\mathcal{I} \times W; M_{(0,\pm k^*)}),$$

Define the space $\mathcal{L} := \bigoplus_{j \in \mathbb{Z}} \mathscr{L}^2(\mathbb{R}; \mathcal{C}^n)$, with norm $||w||_{\mathscr{L}}^2 = \sum_{j \in \mathbb{Z}} ||w_j||_{\mathscr{L}^2(\mathbb{R})}^2$. Define the operator $D : \mathcal{I} \times M_{(0,\pm k^*)} \to \mathcal{L}$ by $D(\epsilon, (v_j)) = (\mathscr{L}_j^{\bar{d}}(\epsilon, v_j))_{j \in \mathbb{Z}^*}$. Notice that we can write the mapping T as

$$T(\epsilon, v) = D(\epsilon)^{-1}[S(\mathcal{N}[\epsilon, S^{-1}(v)])] \tag{C.2}$$

where $D(\epsilon)^{-1}$ is the inverse with respect to the space $M_{(0,\pm k^*)}$, i.e., $D(\epsilon)^{-1}(\cdot) = (\mathscr{L}_j^{-1}(\epsilon)(\cdot))_{j \in \mathbb{Z}^*}$. The result will be a consequence of the following lemma:

**Lemma 7.** *$D \in \mathscr{C}^1(\mathcal{I} \times V); \mathcal{L})$, where $V \subset M_{(0,\pm k^*)}$ is an open subset containing the origin.*

*Proof.* . Fix $(\tilde{\epsilon}, (u_j)) \in \mathcal{I} \times V$, for $V \subset M_{(0,\pm k^*)}$ small neighborhood. For $(\epsilon, (v_j)) \in \mathcal{I} \times M_{(0,\pm k^*)}$, we estimate $D(\epsilon + \tilde{\epsilon}, u + v) - D(\tilde{\epsilon}, u) - \epsilon \partial_\epsilon D(\tilde{\epsilon}, u) - D(\tilde{\epsilon}, v)$, where $\partial_\epsilon D(\tilde{\epsilon}, u) = (\partial_\epsilon \mathscr{L}_j^{\bar{d}}(\epsilon, u))_{j \in \mathbb{Z}^*}$. As the higher order coefficients of the operator $\mathscr{L}^{\bar{d}}$ do not depend on $\epsilon$ we can use the linearity to write the above expression as

$$\left(\partial_x([A(\epsilon + \tilde{\epsilon}) - A(\tilde{\epsilon}) - \epsilon \partial_\epsilon A(\tilde{\epsilon})] u_j)\right)_{j \in \mathbb{Z}} + i\left(j[B(\epsilon + \tilde{\epsilon}) - B(\tilde{\epsilon}) - \epsilon \partial_\epsilon B(\tilde{\epsilon})] u_j\right)_{j \in \mathbb{Z}} +$$
$$+ \left(\partial_x([A(\epsilon + \tilde{\epsilon}) - A(\tilde{\epsilon})] w_j)\right)_{j \in \mathbb{Z}} + i\left(j[B(\epsilon + \tilde{\epsilon}) - B(\tilde{\epsilon})] w_j\right)_{j \in \mathbb{Z}}$$

Notice that $||(j u_j)||_{\mathcal{L}} \leq ||u||_M$ and $||(\partial_x u_j)||_{\mathcal{L}} \leq ||u||_M$ and. In which concerns to the coefficients, the uniform continuity of over compact sets, the boundedness of $\bar{u}$ and $\bar{u}'$ asserts that all the following functions - $||[A(\epsilon+\tilde{\epsilon})-A(\tilde{\epsilon})-\epsilon\partial_\epsilon A(\tilde{\epsilon})]||_{\mathscr{L}^\infty(\Omega)}$, $||B(\epsilon+\tilde{\epsilon})-B(\tilde{\epsilon})-\epsilon\partial_\epsilon B(\tilde{\epsilon})||_{\mathscr{L}^\infty(\Omega)}$, $||\partial[A(\epsilon+\tilde{\epsilon})-A(\tilde{\epsilon})-\epsilon\partial_\epsilon A(\tilde{\epsilon})]||_{\mathscr{L}^\infty(\Omega)}$ and $||[B(\epsilon+\tilde{\epsilon})-B(\tilde{\epsilon})]||_{\mathscr{L}^\infty(\Omega)}$ - are $o(|\epsilon|)$ uniformly in $j \in \mathbb{Z}^*$. Therefore, $||D(\epsilon + \tilde{\epsilon}, u + v) - D(\tilde{\epsilon}, u) - \epsilon \partial_\epsilon D(\tilde{\epsilon}, u) - D(\tilde{\epsilon}, v)||_{\mathcal{L}} = o(|\epsilon| + ||v||_M)$.
$\square$

*Proof of theorem C.1.* Applying the implicit function theorem in Banach spaces [Ham82, section I.5] we obtain the $(\epsilon, v) \mapsto D(\epsilon)^{-1}(v) \in \mathscr{C}^1(\mathcal{I} \times \mathcal{L}; M_{(0,\pm k^*)})$. As $Range(S) \subset \mathcal{L}$, the result follows by tracking differentiability of the composed mappings in the equation (C.2). $\square$

## Appendix D. Fredholm Alternative via Dunford integrals

We outline a simple and general way of getting a quantitative Fredholm alternative for unbounded operators via Dunford integrals and the Calculus of Residues. Our object is to establish the following Fredholm Alternative for general unbounded operators:

**Theorem D.1.** *Let L be a closed, densely defined, operator on a Banach space B, with associated norm $|\cdot|$, and let $\lambda = 0$ be an isolated finite-multiplicity eigenvalue of L, with associated total eigenprojection $\Pi$. Then, $Lu = f$ is soluble if and only if $\Pi f$ lies in the range of $\Pi L = L\Pi L$, in which case there is a solution with $|u| \leq C|f|$. In particular, $Lu = f$ is soluble if $\Pi f = 0$, with a unique solution $|u| \leq C|f|$ in $Range(Id - \Pi)$, that is, L has a bounded inverse on $Range(Id - \Pi)$.*

**Remark 13.** *This may be recognized as the usual Fredholm Alternative in a direct form avoiding reference to dual spaces. When $\lambda = 0$ is semisimple, the solvability criterion reduces to $\Pi f = 0$.*

The issue here is that the resolvent $(\lambda - L)^{-1}$ is not necessarily compact for any value of $\lambda \in \mathbb{C}$, and so the usual Fredholm framework based on compactness arguments is not available. We substitute for this a simple Dunford integral computation. Recall the standard resolvent identity

$$(\lambda - L)^{-1} = (L/\lambda)(\lambda - L)^{-1} + Id/\lambda, \tag{D.2}$$

and the definition by Dunford integrals of the complementary eigenprojections

$$\Pi := \frac{1}{2\pi i} \oint_{\partial B(0,r)} (\lambda - L)^{-1} d\lambda, \tag{D.3}$$

$$\tilde{\Pi} := Id - \Pi = Id - \frac{1}{2\pi i} \oint_{\partial B(0,r)} (\lambda - L)^{-1} d\lambda. \tag{D.4}$$

Define for $r > 0$ sufficiently small,

$$L^\dagger := -\frac{1}{2\pi i} \oint_{\partial B(0,r)} \lambda^{-1} (\lambda - L)^{-1} d\lambda. \tag{D.5}$$

**Proposition 18.** $L^\dagger$ *is well-defined from $B \to D(L)$, bounded from $B \to B$, with $LL^\dagger = L^\dagger L = \tilde{\Pi}$.*

*Proof.* Boundedness and invariance with respect to $r > 0$ sufficiently small follow by analyticity of the resolvent $(\lambda - L)^{-1}$ on the resolvent set, and the assumption that $\lambda = 0$ is an isolated eigenvalue. Applying $L$ to (D.5) and using commutation of $L$ and its resolvent and (D.2), we obtain, formally,

$$L^\dagger L = LL^\dagger = -\frac{1}{2\pi i} \oint_{\partial B(0,r)} (L/\lambda)(\lambda - L)^{-1} d\lambda$$

$$= -\frac{1}{2\pi i} \oint_{\partial B(0,r)} (\lambda - L)^{-1} d\lambda + \frac{1}{2\pi i} \oint_{\partial B(0,r)} (1/\lambda) d\lambda$$

$$= -\Pi + Id = \tilde{\Pi},$$

To make this computation rigorous, we may apply $L$ to a convergent series (with respect to $|\cdot|$) of approximating Riemann sums for $L^\dagger f$, $f \in B$, to obtain a convergent series of approximating Riemann sums for $\tilde{\Pi} f$. Appealing to closedness of $L$, we obtain $L^\dagger f \in D(L)$, with $LL^\dagger f = \tilde{\Pi} f$. □

*Proof of Theorem D.1.* For $\Pi f = 0$, we have $Lu = \tilde{\Pi} f = f$ for $u = L^\dagger f$, with $|u| \leq |L^\dagger||f| \leq C|f|$. For $\tilde{\Pi} f = 0$, we have $Lu = f$ if and only if $\Pi Lu = \Pi f = f$, or if and only if $\Pi f$ is in the range of $\Pi L$, in which case we may again (this time appealing to finite-dimensional theory) find a solution satisfying a uniform bound $|u| \leq C|f|$. Combining these facts, we obtain the full result by linear superposition. Finally, $L^\dagger L = LL^\dagger = \tilde{\Pi}$ implies that $L^\dagger$ is a bounded inverse of $L$ on $Range \tilde{\Pi}$. □

**Remark 14.** *Examining the above arguments, we see that it is the property of closedness of $L$ that substitutes here for the compactness of the resolvent used in standard Fredholm theory.*

**Acknoledgements.** This work was done as part of the thesis of Rafael Monteiro under the supervision of prof. Kevin Zumbrun, whom Rafael would like to thank for the fruitful conversations and help (especially with appendix D). He would also like to thanks Miguel Rodrigues for helpful comments.

## References

[AS12] Franz Achleitner and Peter Szmolyan. Saddle-node bifurcation of viscous profiles. Phys. D, 241(20):1703–1717, 2012.

[CH82] Shui Nee Chow and Jack K. Hale. Methods of bifurcation theory, volume 251 of Grundlehren der Mathematischen Wissenschaften [Fundamental Principles of Mathematical Science]. Springer-Verlag, New York, 1982.

[Cha61] S. Chandrasekhar. Hydrodynamic and hydromagnetic stability. The International Series of Monographs on Physics. Clarendon Press, Oxford, 1961.

[CR71] Michael G. Crandall and Paul H. Rabinowitz. Bifurcation from simple eigenvalues. J. Functional Analysis, 8:321–340, 1971.

[CR73] Michael G. Crandall and Paul H. Rabinowitz. Bifurcation, perturbation of simple eigenvalues and linearized stability. Arch. Rational Mech. Anal., 52:161–180, 1973.

[Cra77] M. G. Crandall. An introduction to constructive aspects of bifurcation and the implicit function theorem. In Applications of bifurcation theory (Proc. Advanced Sem., Univ. Wisconsin, Madison, Wis., 1976), pages 1–35. Publ. Math. Res. Center, No. 38. Academic Press, New York, 1977.

[FT08] Heinrich Freistühler and Yuri Trakhinin. On the viscous and inviscid stability of magnetohydrodynamic shock waves. Phys. D, 237(23):3030–3037, 2008.

[Ham82] Richard S. Hamilton. The inverse function theorem of Nash and Moser. Bull. Amer. Math. Soc. (N.S.), 7(1):65–222, 1982.

[Hen81] Daniel Henry. Geometric theory of semilinear parabolic equations, volume 840 of Lecture Notes in Mathematics. Springer-Verlag, Berlin, 1981.

[Kno96] Edgar Knobloch. Symmetry and instability in rotating hydrodynamic and magnetohydrodynamic flows. Phys. Fluids, 8(6):1446–1454, 1996.

[MZ05] Guy Métivier and Kevin Zumbrun. Large viscous boundary layers for noncharacteristic nonlinear hyperbolic problems. Mem. Amer. Math. Soc., 175(826):vi+107, 2005.

[MZ09] Guy Métivier and Kevin Zumbrun. Existence of semilinear relaxation shocks. J. Math. Pures Appl. (9), 92(3):209–231, 2009.

[PYZ03] A. Pogan, J. Yao, and K. Zumbrun. O(2) hopf bifurcation of viscous shock waves in a channel. preprint, 2003.

[Sat76] D. H. Sattinger. On the stability of waves of nonlinear parabolic systems. Advances in Math., 22(3):312–355, 1976.

[SZ99] Denis Serre and Kevin Zumbrun. Viscous and inviscid stability of multidimensional planar shock fronts. Indiana Univ. Math. J., 48(3):937–992, 1999.

[TL86] Angus E. Taylor and David C. Lay. Introduction to functional analysis. Robert E. Krieger Publishing Co. Inc., Melbourne, FL, second edition, 1986.

[TZ05] Benjamin Texier and Kevin Zumbrun. Relative Poincaré-Hopf bifurcation and galloping instability of traveling waves. Methods Appl. Anal., 12(4):349–380, 2005.

[TZ08a] Benjamin Texier and Kevin Zumbrun. Galloping instability of viscous shock waves. Phys. D, 237(10-12):1553–1601, 2008.

[TZ08b] Benjamin Texier and Kevin Zumbrun. Hopf bifurcation of viscous shock waves in compressible gas dynamics and MHD. Arch. Ration. Mech. Anal., 190(1):107–140, 2008.

[Yos95] Kōsaku Yosida. Functional analysis. Classics in Mathematics. Springer-Verlag, Berlin, 1995. Reprint of the sixth (1980) edition.

Indiana University, Bloomington, IN 47405
*E-mail address*: `radmonte@indiana.edu`